\documentclass[12pt,reqno]{amsart}

\newcommand{\mysection}[1]{\section{#1}
      \setcounter{equation}{0}}

\usepackage{color}

\newtheorem{theorem}{Theorem}[section]
\newtheorem{lemma}[theorem]{Lemma}

\newtheorem{corollary}[theorem]{Corollary}

\theoremstyle{definition}
\newtheorem{assumption}{Assumption}[section]

\newtheorem{example}{Example}[section]

\theoremstyle{remark}
\newtheorem{remark}{Remark}[section]

\newcommand\bR{\mathbb{R}}

\newcommand\frK{\mathfrak{K}}

\newcommand\cK{\mathcal{K}}
\newcommand\cL{\mathcal{L}}

\newcommand\cQ{\mathcal{Q}}

\makeatletter
 \newcommand{\sumstar}%
 {\operatornamewithlimits{\sum@\kern-.2em\raise1ex\hbox{*}}}
 \makeatother

\begin{document}

\title[First derivatives estimates]
{First derivatives estimates for finite-difference schemes}

\author[I. Gy\"ongy]{Istv\'an Gy\"ongy}
\address{School of Mathematics,
University of Edinburgh,
King's  Buildings,
Edinburgh, EH9 3JZ, United Kingdom}
\email{gyongy@maths.ed.ac.uk}

\author[N.  Krylov]{Nicolai Krylov}%
\thanks{The work of the second author was partially supported
by NSF grant DMS-0653121}
\address{127 Vincent Hall, University of Minnesota,
Minneapolis,
       MN, 55455, USA}
\email{krylov@math.umn.edu}

\subjclass{65M06,39A70}
\keywords{Cauchy problem,  
  finite differences, first derivatives estimates}

\begin{abstract}We give sufficient conditions
under which solutions of discretized in space
second-order
parabolic and elliptic equations, perhaps degenerate,
admit estimates of the first derivatives in the space
variables independent of the mesh size.
\end{abstract}

\maketitle

\mysection{Introduction} 
                  \label{section02.04.06}

This is the first article out of a series of two
devoted to estimating space derivatives of
solutions of discretized in space second-order
parabolic and elliptic equations. We allow
equations to degenerate and to become just first-order
equations. In the present article we only deal with the 
first-order derivatives. In the second part of this project we
will prove higher-order derivatives estimates
and apply them to showing a method
of accelerating finite-difference approximations
 to any given rate  for equations 
in the whole space. 

 Numerical approximations for
linear and quasilinear
partial differential equations is a rather old and well developed
area. We refer the reader to \cite{DK} and the references
therein, following which
one can track down original papers by 
D. Aronson, L. Bers, R. Courant-K. 
 Friedrichs-H. Lewy,
J. Douglas, F. John,
O. Ladyzhenskaya, P. Lax,  H. Levy, L. Liusternik,
I. Petrovskii,  
and many many others to which we only add 
\cite{Li} and two more papers
\cite{Gy} and \cite{Yo} where
discrete methods are applied to approximate stochastic
partial differential equations.

A major difference
of this article from 
all above mentioned ones is that we focus on
investigating the smoothness of approximating solutions
rather than on convergence only.
For each point $x\in\bR^{d}$ we move the original
grid in such a way that $x$ becomes a grid point.
This allows us to define the approximate solution
in all of $\bR^{d}$ rather than only on the grid 
and we investigate how smooth the approximating solution
is with respect to $x$. We estimate true derivatives rather
than their difference approximations.

Estimating the sup norms of the
first-order derivatives for solutions
of finite-difference
 schemes for linear and {\em fully\/} nonlinear
second-order  
{\em degenerate\/} 
equations plays a major role
in estimating the rate of convergence of
approximating solutions to the true solution
in the sup norm.
The most general results for fully nonlinear equations
concerning the rates can be found in \cite{DK1}
and in the references therein. A recent development 
in the issue of estimating the Lipschitz constant 
 and second-order differences for
approximating solutions for fully nonlinear equations
without applications to estimating
the rate of convergence is presented in \cite{Kr07}.
Before, the Lipschitz constants and higher order 
derivatives
estimates were obtained in \cite{DK} 
 for time-space discretization of 
linear degenerate  parabolic equations
 and in \cite{DK2} for fully nonlinear 
equations. They are also
applied to estimating the rate of convergence.
In a sense the present article is close to \cite{DK}.
However, here we only deal with the first-order 
derivatives
estimates and for equations discretized only 
in the space
variable. We introduce a new type of sufficient 
conditions
for obtaining the estimates
(see Assumption \ref{assumption 11.22.11.06}
and the discussion in Section
\ref{section 4.9.1}). These conditions are
much weaker 
and more detailed than the corresponding ones in 
\cite{DK}.
Our method is also somewhat  different. Instead
of considering just the sum of squares 
of the difference
increments along the mesh we add to it the 
square of the full gradient with a small 
constant factor. 
This allows us to estimate 
the gradient.

In this connection it is worth noting that
such an estimate for  
finite-differences approximations of the first-order
directional derivatives
in $x$ is claimed 
in Theorem 4.1 of \cite{DK} under some conditions,
which are always satisfied if the equation
is uniformly nondegenerate even if $c$ (see
\eqref{07.9.25.03}) is small. However, in this case,
actually, the result of Theorem 4.1 of \cite{DK}
is only proved for the derivatives along the mesh. 
This is rather
harmless if the vectors on the mesh span the whole space,
but excludes the cases when the mesh lies in a subspace, which 
happens, for instance, if
  we are dealing with, say uniformly nondegenerate 
equations whose coefficients depend
  on a parameter and we want to estimate the 
the finite-differences of their
solutions with respect to the parameter 
by considering it as just another space variable.
In that case no second-order derivatives
with respect to the parameter enters the limit equation,
 the assumption that it is uniformly nondegenerate
with respect to the original space variables does not help,
and we need to have $c $ be large in order to 
rely on Theorem 4.1 of \cite{DK}. Our results are free
from this flaw, see Remarks \ref{remark 07.11.25.1}
and \ref{remark 3.8.1}.

To understand faster the method and the results 
of the article, we advise the reader concentrate
only on the parabolic case and
assume
that the limit equation is uniformly nondegenerate.
Then apart from the Lipschitz continuity
nothing  else (see   Remarks \ref{remark 07.9.18.8}
and \ref{remark 11.10.06})
 is required
for Theorem \ref{theorem 4.7.2} to hold.

For the general equations
our conditions (see 
Remarks \ref{remark 07.9.18.5} and \ref{remark 07.9.18.9})
capture the main features
of the corresponding conditions known from the theory
of PDE. Namely, roughly speaking, we need the first-order
derivatives
in any direction along the mesh
of the coefficients to be dominated either by the diffusion
coefficients along the same direction or
by $c$ or else by the drift term if it is sufficiently
``monotone".
It is worth noting that along
 the way we discover the necessity
of the diffusion coefficients to have a special form
and the usefulness of adding a diffusion term
with a coefficient proportional to the mesh step
into {\em approximating\/} equation. This reminds
the method of artificial diffusion, although,
as far as we understand,
the artificial diffusion is
usually added to the
{\em original\/} differential equation.

Of course, in the same way as in \cite{DK},
the results of the present article
lead to the rate of convergence of order $h^{1/2}$
of approximating solutions to the true solution.
However, for brevity we do not say more about this
issue only adding that in general
our finite-difference equations
need not be related in any way to 
a partial differential equation.

Our main results are collected in Section 
 \ref{section 13.19.10.07},
which also contains the proofs
of  all of them but Theorem \ref{theorem 4.7.2}, 
which is proved in Section \ref{section 07.9.25.1}.
In Section  \ref{section 13.19.10.07} 
we also state, in a special case,
without proof one of 
the main results of the continuation of the present paper.
Section \ref{section 07.9.25.6} contains a
discussion of our assumptions concerning
the structure of the finite-difference
equations under consideration. The point is
that our equations do not contain mixed
second-order differences and in Section \ref{section 07.9.25.6}
we explain that this is ``almost" the most interesting case.
The final Sections
\ref{section 4.9.1} and \ref{section 07.11.26.1}
are devoted to a rather long and
 detailed discussion of the somewhat
formally stated
Assumption \ref{assumption 11.22.11.06} and showing
that it is natural in many   cases 
alluded to above.

\mysection{Formulation of the main results}
                                 \label{section 13.19.10.07}

We take some numbers  $h, T\in(0,\infty)$ and
 in a cylindrical domain consider the integral equation 
\begin{equation}                      \label{equation}
u(t,x)=g (x)+\int_0^t\big(L u(s,x)+f (s,x)\big)\,ds
\end{equation}
for $u$, where  $ g (x)$ and $ 
f (s,x)$  are
given  real-valued Borel functions of $x
=(x_{1},...,x_{d})\in\bR^d$ and  
$(s,x)\in H_T:=[0,T]\times\bR^d$,  
respectively, and 
 $L $ is a linear operator 
given by 
$$
L \varphi(t,x)=L_{h}\varphi(t,x)=L ^0\varphi(t,x)
-c(t,x)\varphi(x), 
$$
\begin{equation}
                                     \label{07.9.25.03}
L ^0\varphi(t,x)=L_{h}^0\varphi(t,x)=\frac{1}{h }
\sum_{\lambda\in\Lambda_{1}}
q_{\lambda}(t,x )\delta_{\lambda}\varphi(x ) 
+ \sum_{\lambda\in\Lambda_{1}} 
p_{\lambda}(t,x )\delta_{\lambda}\varphi(x ), 
\end{equation}
for functions $\varphi$ on $\bR^d$. Here 
$\Lambda_{1}$ 
is a finite  subset  of $\bR^d$ such that
$0\not\in\Lambda_{1}$, 
  $p_{\lambda}(t,x )$, $q_{\lambda}(t,x )$ are
real-valued functions  of $(t,x)\in H_T$
given for each 
$\lambda\in\Lambda_{1}$, and
$$
\delta_{\lambda}\varphi(x )=\delta_{h,\lambda}\varphi(x )=
\frac{1}{h}(\varphi(x +h\lambda)-\varphi(x )),
\quad\lambda\in\Lambda_{1}.
$$ 

Let $m\geq0$ be an integer and let $K_{1}\in[1,\infty)$ be
a constant. 
Introduce
$$
\chi_{\lambda}=\chi_{h,\lambda}:=q_{\lambda} +hp_{\lambda}. 
$$
We make the following assumptions. 

\begin{assumption}                \label{assumption 16.12.07.06}
The 
functions $p$, $q$, $c$,  
$f$, and $g$ and their
derivatives in $x$   up to order $m$ are  bounded on $H_T$ and
continuous in $x$.
\end{assumption}

\begin{assumption}  
                            \label{assumption 1.26.11.06}   
For all $(t,x)\in H_T$
and $\lambda\in\Lambda_{1}$,
$$
\chi_{\lambda}(t,x )\geq 0.
$$
There exists a constant $c_{0}>0$ such that $c\geq c_{0}$.
\end{assumption}  
\begin{remark}
                                     \label{remark 07.9.18.8}
The above assumption: $c\geq c_{0}>0$, is almost irrelevant
if we only consider \eqref{equation} on a finite time interval.
Indeed,  if $c$ is just bounded, say $|c|\leq C=\text{const}$, by 
introducing a new function $v(t,x)=u(t,x)e^{-2Ct}$
we will have an equation for  $v$ similar to 
\eqref{equation} with $L^{0}v-(c+2C)v$ and $fe^{-2Ct}$
in place of $Lu$ and $f$, respectively. Now for the new $c$ we have
  $c+2C\geq C$.
\end{remark}

\begin{remark}
                                     \label{remark 3.8.001}
  
Introduce the following 
symmetry condition:\medskip

\noindent(S) We have $\Lambda_{1}=-\Lambda_{1}$ and $q_{\lambda}=
q_{-\lambda}$ on $\Lambda_{1}$.\medskip

Obviously under condition (S) we have
$$
h^{-1}\sum_{\lambda\in\Lambda_{1}}q_{\lambda}(t,x)
\delta_{\lambda}\varphi(x)
=(1/2)\sum_{\lambda\in\Lambda_{1}}q_{\lambda}(t,x)
\Delta_{\lambda}\varphi(x),
$$
where
$$
\Delta_{\lambda}\varphi(x)=h^{-2}
(\varphi(x+h\lambda)-2\varphi(x)+\varphi(x-h\lambda)). 
$$
\end{remark}

Take a   function
$\tau_{\lambda}$ defined on $\Lambda_{1}$ 
 taking values in $[0,1]$,  
and for  $\lambda\in\Lambda_{1}$ 
  introduce
the   operators
$$
T_{ \lambda}\varphi(x)=T_{ h,\lambda}\varphi(x)=
\varphi(x+h\lambda),\quad
\bar{\delta}_{\lambda}=\tau_{\lambda}h^{-1}
(T_{ \lambda}-1) . 
     $$ 
It is worth noticing that in most applications we take
$\tau_{\lambda}\equiv1$ on $\Lambda_{1}$.
 However, there are cases
(see Remark \ref{remark 07.11.25.1}) 
in which it is useful to have
some flexibility in changing $\tau_{\lambda}$.
 
For uniformity of notation we also introduce
$\Lambda_{2}$ as the set of fixed
 distinct vectors $\ell^{1},...,
\ell^{d}$ none of which is in $ 
\Lambda_{1}$ and   define
$$
\bar{\delta}_{ \ell^{i}}=\bar{\delta}_{h,\ell^{i}}
=\tau_{0}D_{i}:=\tau_{0}
\partial/\partial x_{i},\quad T_{ \ell^{i}}=1,\quad
\Lambda=\Lambda_{1}\cup\Lambda_{2},
$$
where $\tau_{0}\in[0,1]$ 
is a fixed constant.
Observe that we allow $\tau$ to be zero
in order to cover some results from \cite{DK}. 

For $\mu\in\Lambda $   we set
$$
Q \varphi =h^{-1}\sum_{\lambda\in\Lambda_{1}}
q_{\lambda} \delta_{ \lambda}
\varphi,
\quad
Q_{  \mu}\varphi =h^{-1}\sum_{\lambda\in\Lambda_{1}}
  (\bar{\delta}_{ \mu}q_{\lambda})\delta_{ \lambda}
\varphi  ,
$$
$$
P \varphi = \sum_{\lambda\in\Lambda_{1}}
p_{\lambda} \delta_{ \lambda}
\varphi,
\quad
P_{ \mu}\varphi = \sum_{\lambda\in\Lambda_{1}}
 (\bar{\delta}_{ \mu}p_{\lambda} )\delta_{ \lambda}
\varphi,
$$
$$
L^{0}_{ \mu}
=Q_{\mu }
+P_{\mu}.
$$

Below $B(\bR^{d})
$ is the set of  bounded Borel functions on
$\bR^{d}$ and $\frK$ is the set of bounded operators
$\cK =\cK (t)$ mapping $B(\bR^{d})
$ into itself preserving  the cone
of nonnegative functions and satisfying $\cK 1\leq1$. 
 We will often make use of the  
simple fact that for any 
$\cK_1,\cK_2\in\frK$ and nonnegative 
functions $\alpha$, $\beta$ on $\bR^d$,  
$$
\alpha\cK_1+\beta\cK_2=(\alpha+\beta)\cK_3
$$ 
with 
$$
\cK_3:=\tfrac{\alpha}{\alpha+\beta}\cK_1
+\tfrac{\beta}{\alpha+\beta}\cK_2
\in\frK, \quad \big(\tfrac{0}{0}:=0\big).
$$

\begin{assumption}         \label{assumption 11.22.11.06} 
We have $m\geq1$ and there
exist  a  constant    
$\delta \in(0,1 ] $ 
 and an operator $\cK=\cK_{h}\in\frK$,
 such that  
\begin{equation}
                                           \label{3.24.1}
2  \sum_{\lambda\in\Lambda }
(\bar{\delta}_{ \lambda}\varphi)
L^{0}_{ \lambda}T_{ \lambda}\varphi\leq
 \sum_{\lambda \in\Lambda}
\cQ( \bar{\delta}_{ \lambda}
\varphi)
+K_{1}\cQ(\varphi)
+
2(1-\delta)c\cK \big(\sum_{\lambda\in\Lambda} 
|\bar{\delta}_{\lambda}\varphi|^{2}\big)  
\end{equation}
on $H_{T}$
for all smooth functions $\varphi$, where
$$
\cQ(\varphi)=\sum_{\mu\in\Lambda_{1}}\chi_{\mu}
| \delta_{\mu}\varphi|^{2}.
$$
\end{assumption}
 
It is worth noting that Assumption
\ref{assumption 11.22.11.06} is automatically      
satisfied if  $q_{\lambda}$ and $p_{\lambda}$    
         are independent of $x$. 
There are a few more cases when it is satisfied as well.
We discuss some of them 
here and  in Section \ref{section 4.9.1} 
only mentioning  
right away three situations.  

\begin{remark}                    
                                    \label{remark 11.10.06}
Let Assumptions \ref{assumption 16.12.07.06} 
and \ref{assumption 1.26.11.06} 
hold with $m\geq1$. Assume that 
$\Lambda_1=-\Lambda_1$, $Dq_{-\lambda}=Dq_{\lambda}$ 
and $q_{\lambda}\geq\kappa$
for all $\lambda\in\Lambda_1$, where $\kappa>0$ 
is some constant.  
Then Assumption
\ref{assumption 11.22.11.06} is satisfied as well 
with  $\delta$ as close to $1$ as we wish, 
with $\tau_{\lambda}\equiv1$ on $\Lambda_{1}$,
appropriate $\tau_{0}>0,K_{1}$,
unit $\cK$, and all small $h$. 
\end{remark}
We will prove this remark at the end of this section.
In Remark \ref{remark 07.9.18.7} we show that if we
have $m\geq2$
 and the symmetry condition
(S) is satisfied, then  in the above remark  
the condition $\kappa>0$
can be replaced with $\kappa=0$, provided that $c_{0}$
is large enough 
(this time we need not assume that $h$ is small).
In Remark \ref{remark 07.9.18.7} we also show that
  the condition $m\geq2$ 
can be replaced with $m=1$ provided that 
  $\sqrt{q_{\lambda}}$ are Lipschitz
continuous in $x$ with a constant independent of $t$.
  In that case 
again Assumption
\ref{assumption 11.22.11.06} is satisfied  
for $\delta=1/10$ and appropriate $K_{1},
 \tau_{0}>0 $, 
provided that $c_{0}$
is large enough. 

 As we have seen in Remark \ref{remark 07.9.18.8},
the condition that $c_{0}$ be large is, actually,
harmless as long as we are concerned with equations
on a finite time interval.

Fix a domain $Q\subset\bR^{d}$ and introduce 
$$
Q^{o}
=\{
x\in Q:x+\lambda h\in Q\quad\forall\lambda\in\Lambda_{1}
\},
\quad\delta Q=Q\setminus Q^{o},
$$
$$
Q_{T}=[0,T]\times Q,\quad Q^{0}_{T}=
[0,T]\times Q^{0},\quad\delta_{x}Q_{T}=
[0,T]\times\delta Q,
$$
$$
\delta'Q_{T}=(\{0\}\times Q)\cup\delta_{x}Q_{T}.
$$

Our first main result is formulated as follows,
where by $Du$ we mean the gradient of $u$ with respect
to $x$.  Observe that
the 
main case that $Q=\bR^{d}$ 
is not excluded and in this case
assumption (ii) of Theorem \ref{theorem 4.7.2} below
can be checked on the basis of  
Theorem \ref{theorem 07.9.12.1}.
A typical and the most reasonable
 application of Theorem
\ref{theorem 4.7.2} when $Q$ is a proper domain
is the case that $\tau_{0}=0$.

\begin{theorem}
                                   \label{theorem 4.7.2}
(i) Let
Assumptions \ref{assumption 16.12.07.06} through 
\ref{assumption 11.22.11.06}    be satisfied and
let $u$ be a bounded function on 
$H_{T}$ satisfying
\eqref{equation} in $  Q _{T}$.

(ii) Assume that $u$  and $Du$ are 
bounded and continuous
in $Q_{T}$.

Then
in $Q_{T}$ we have 
\begin{equation}
                                             \label{4.8.03}
|u|+\tau_{0}|D u|+U \leq N  
(  F_{1}+\sup_{\delta'Q_{T}}(|u|+\tau_{0}|D u|+U )),  
\end{equation}
where
$$
U=\big(\sum_{\lambda\in\Lambda_{1}}|\bar{\delta}_{\lambda}
u|^{2}\big)^{1/2},\quad
F_{1}=\sup_{H_{T}}(|f|+
|D f|),
$$
  and $N $
depends only on $\delta,c_{0},  K_{1} $, 
$\sup_{H_{T}}|Dc|$, and
$$
|\Lambda_{1}|^{2}:=\sum_{\lambda\in\Lambda_{1}}|
\tau_{\lambda}\lambda|^{2}.
$$ 
\end{theorem}

\begin{remark}
                                   \label{remark 07.11.26.5}
We will see from the proof that, if $\tau_{0}=0$,
then Theorem \ref{theorem 4.7.2} holds without the
assumption that $Du$ exists let alone continuous. 

\end{remark}

In case $Q=\bR^{d}$  
assumption (ii) of the following result
is often satisfied due to Theorem \ref{theorem 07.9.12.1}.

It is worth noting that if  $Q=\bR^{d}$, then $\delta Q=
\emptyset$ and for any function $\varphi$ we set
$$
\sup_{\emptyset}\varphi:=0.
$$

\begin{theorem}
                                   \label{theorem 07.9.25.2}
(i) Let  
Assumptions  \ref{assumption 16.12.07.06} 
through \ref{assumption 11.22.11.06}
 be satisfied. Suppose that $q_{\lambda}$, $p_{\lambda}$,
$c$, and $f$ are independent of $t$.  

(ii) Assume that in $\bR^{d}$
there exists a 
bounded
function $u=u(x)$ which is bounded and continuous
in $Q$ along with $Du$ and such that
$$
Lu+f=0\quad\text{in}\quad Q.
$$
Then in $Q$ we have 
$$
|u|+\tau_{0}|Du|+U\leq N(F_{1}
+ \sup_{\delta Q}(|u|+\tau_{0}|Du|
+U),
$$
where $U$ is the same as in Theorem \ref{theorem 4.7.2},
$$
F_{1}=\sup_{\bR^{d}}(|f|+
|D f|),
$$
  and $N $
depends only on $\delta,c_{0},  K_{1} $, 
$\sup_{\bR^{d}}|Dc|$, and
 $|\Lambda_{1}|$.
\end{theorem}

Proof. Take $\nu=c_{0}/2$, so that $c-\nu\geq c_{0}/2$,
and observe that in $Q_{T}$ the function 
$v(t,x):=u(x)e^{\nu t} $
satisfies
\begin{equation}
                                     \label{07.9.25.10}
\frac{\partial}{\partial t}v=L^{0}v-(c-\nu)v+e^{\nu t}f.
\end{equation}
By Theorem \ref{theorem 4.7.2} for $x\in Q$ 
and obvious meaning of $V$ we have
$$
e^{\nu T}(|u(x)|+\tau_{0}|Du(x)|+U(x))=
|v(T,x)|+\tau_{0}|Dv(T,x)|+V(T,x)
$$
$$
\leq
N e^{\nu T}[F_{1}+\sup_{\delta Q}(|u|+\tau_{0}|Du|+U)]
+N\sup_{Q}(|v |+\tau_{0}|Dv |+V)(0,y).
$$
By multiplying the extreme terms by $e^{-\nu T}$
and letting $T\to\infty$, we get the result. The theorem is
proved.
 
\begin{remark}It is worth noticing that in the above
theorems it suffices that 
\eqref{3.24.1} be satisfied
only in  $Q^{o}_T$. 
\end{remark}

\begin{theorem}
                                   \label{theorem 07.9.12.1}

(i) Let
Assumption  \ref{assumption 16.12.07.06}  be satisfied.
Then there exists a unique bounded solution $u$ of 
\eqref{equation} in $H_{T}$. Moreover,  
all derivatives  in $x$ of $u$ of order $\leq m$
are bounded and continuous in $H_{T}$. 

(ii)  Let
Assumptions  \ref{assumption 16.12.07.06} 
through \ref{assumption 11.22.11.06}
 be satisfied. Suppose that $q_{\lambda}$, $p_{\lambda}$,
$c$, and $f$ are independent of $t$. 
Then there exists a unique bounded solution $u=u(x)$ of the
equation
\begin{equation}
                                       \label{07.9.25.8}
Lu+f=0\quad\text{in}\quad\bR^{d}.
\end{equation}
  Moreover,  $u$ and  $Du$ are bounded and continuous
in $\bR^{d}$. 
\end{theorem}

Proof. (i) Let $C^{m}$ be the space of functions
on $\bR^{d}$ which are bounded and continuous along with
all derivatives up to order $m$. We endow $C^{m}$ with an 
appropriate sup norm and in so obtained Banach space,
denoted again by $C^{m}$, consider
the equation
$$
u(t)=g+\int_{0}^{t}(A(s)u(s)+f(s))\,ds,
$$
where $f(s)=f(s,\cdot)$ and $A(s)$ are operators in $C^{m}$ given by
$$
A(s)\varphi(x)=h^{-1}\sum_{\lambda\in\Lambda_{1}}
\chi_{\lambda}(s,x)\delta_{\lambda}\varphi(x)-c(s,x)\varphi(x).
$$

Owing to Assumption  \ref{assumption 16.12.07.06}
$$
\|A(s)\varphi\|_{C^{m}}\leq N\|\varphi\|_{C^{m}}
$$
with $N$ independent of $s$ and $\varphi$. Hence, our result is a direct
consequence of the general theorem about ordinary differential
equations in Banach spaces (for proving 
uniqueness we take $m=0$).

(ii) By assertion (i) for any $T$ there exists a unique bounded
and continuous in $H_{T}$ 
solution $v(t,x)$ of the problem
$$
\frac{\partial}{\partial t}v(t,x)=
 (L+\nu)v(t,x) \quad t>0,\quad v(0,x)=f(x),
$$
where $\nu=c_{0}/2$. In addition, 
 $Dv$ is bounded and continuous in $H_T$ 
for each $T$.  
By Theorem \ref{theorem 4.7.2},
$v$ and  $Dv$ are bounded
and continuous in $H_{\infty}$. Define 
$$
u(x)=\int_{0}^{\infty}e^{-\nu t}v(t,x)\,dt.
$$
Then the rules of differentiating under the integral sign
and the dominated convergence theorem
show that $u$ and $Du$ are bounded and continuous.
Furthermore, integrating by parts, we see that
$$
Lu(x)=\int_{0}^{\infty}e^{-\nu t}Lv(t,x)\,dt
=\int_{0}^{\infty}e^{-\nu t}\big[
\frac{\partial}{\partial t}v(t,x)-\nu v(t,x)\big]\,dt
=-f(x),
$$
so that $u$ satisfies \eqref{07.9.25.8}.

To prove uniqueness of bounded solutions of
\eqref{07.9.25.8} we use Lemma \ref{lemma 3.6.1}
which is proved in Section \ref{section 07.9.25.1}.
If $w$ is the difference of two bounded
solutions of \eqref{07.9.25.8},
then 
$v(t,x):=w(x)e^{c_{0}t}$ 
satisfies \eqref{07.9.25.10} with $\nu=c_{0}$ and $f=0$. 
Since $c-\nu\geq  0$,
by Lemma \ref{lemma 3.6.1}  
$$
v(t,x)\leq \sup_{\bR^{d}}v_{+}(0,y)
= \sup_{\bR^{d}}w_{+},\quad
w(x)\leq e^{-c_{0}t} \sup_{\bR^{d}}w_{+}
$$
and by letting $t\to\infty$ we obtain $w\leq0$.
The same inequality holds for $-w$, so that $w=0$,
which proves  uniqueness and finishes 
the proof of the theorem.

\begin{remark}
Let Assumption \ref{assumption 16.12.07.06} hold. 
Then it is easy to see that for the   
bounded solution $u$ of 
\eqref{equation} in $H_{T}$ 
for each fixed 
$h_1>0$ 
$$
\sup_{(t,x)\in H_T}|u(t,x)|\leq N
\sup_{(t,x)\in H_T}|f(t,x)|, 
$$
holds for all $h>h_1$, where $N$ is a constant 
independent of $h$.  
\end{remark}

\begin{remark}
A simple inspection of their proof 
shows that the above theorems remain 
valid if 
$\cK \big(\sum_{\lambda\in\Lambda} 
|\bar{\delta}_{\lambda}\varphi|^{2}\big)$ 
in 
Assumption \ref{assumption 11.22.11.06}
is 
replaced by  
$$
\sup_x\sum_{\lambda\in\Lambda} 
|\bar{\delta}_{\lambda}\varphi|^{2}.
$$  
We do not know how much can be gained 
by such weakening of Assumption \ref{assumption 11.22.11.06}.
On the other hand,
in a subsequent article we will see 
an advantage of using 
operators $\cK\in\frK$. 
\end{remark}

Now we state without proof one of the main results
of the forthcoming paper 
\cite{GK}. As we know
from Remark \ref{remark 11.10.06}
and Theorem \ref{theorem 07.9.12.1} (ii),
under the conditions of Theorem
\ref{theorem 9.4.10.07} (see below),
for each $h>0$,
there exists a unique bounded solution
$u_{h}$
  of 
$$
L_hu+f=0 \quad\text{in $\mathbb R^d$}.
$$
For a fixed integer $k\geq0$ set 
$$
v_h=\sum_{j=0}^{k}b_ju_{2^{-j}h}, 
$$
where 
$$
(b_0,b_1,...,b_k)
:=(1,0,0,...,0)V^{-1}
$$
and $V^{-1}$ is the inverse of the 
Vandermonde matrix with entries
$$
V^{ij}:=2^{-(i-1)(j-1)}, \quad i,j=1,...,k+1.
$$
Consider also the equation
\begin{equation}                    \label{15.05.10.07}
\mathcal Lv+f=0\quad\text{in}\quad \bR^d
\end{equation}
with 
\begin{equation}                    \label{07.10.11.1}
\cL:=a_{ij}D_iD_j+b_iD_i-c, 
\end{equation}
$$
a_{ij}(x):=(1/2)
\sum_{\lambda\in\Lambda_{1}}
q_{\lambda}(x)\lambda_{i}\lambda_{j},
\quad 
b_{i}:=
\sum_{\lambda\in\Lambda_{1}}
p_{\lambda}(x)\lambda_{i}.
$$

One of our main theorems from \cite{GK} 
in a special case reads as follows. 
\begin{theorem}                   \label{theorem 9.4.10.07} 
Let  $m\geq 3(k+1)$ 
 for some integer $k\geq0$.
Let
Assumptions 
\ref{assumption 16.12.07.06}, 
\ref{assumption 1.26.11.06}, 
 and the symmetry assumption (S) be satisfied.
Also assume that $q_{\lambda}(x)\geq\kappa$
for all $\lambda\in\Lambda_1$, $x\in\bR^d$ 
for some constant $\kappa>0$.
         
Then there is a unique bounded solution $v$ to 
\eqref{15.05.10.07} and 
$$
|v_h(x)-v(x)|\leq N h^{k} 
$$
for all $x\in\bR^d$, $h\in(0,h_0]$, and every $h_0>0$, where 
$N$ is a constant depending only on $h_{0}$, 
$k$, $m$,  $\kappa$, $c_0$, $|\Lambda_1|$, and 
on the sup norms of the derivatives of 
$q_{\lambda}$, $p_{\lambda}$, $c$, and 
$f$ 
up to order $m$. 
\end{theorem}

  We obtain this result in 
\cite{GK} by showing that the derivatives  
of $v_h$ in $h$ up to order $k+1$ are bounded 
functions of $h\in(0,h_0]$, which we will prove  
via our estimates on the derivatives of $v_h$ in $x$. The reader 
may wonder why we do not choose the straightforward way 
of estimating the derivatives of $v_h$ in $h$ via 
an `explicit' formula for $v$. 
To test this approach we suggest the reader   try 
to estimate $dv_h/dh$ directly for  
$$
v_h(x)=\tfrac{h^2}{2+h^2}f(x)+
\sum_{n=1}^{\infty}h^2(\tfrac{2}{2+h^2})^{n+1}
Ef(x+\sum_{i=1}^{n}h\varepsilon_i),
$$  
where ${\varepsilon_i}$ are independent random variables 
taking $1$ and $-1$ with probability $1/2$, without 
noticing that $v_h$ is the bounded solution of 
$$
\tfrac{1}{h^2}(u(x+h)-2u(x)+u(x-h))-u(x)+f(x)=0, \quad x\in\bR.
$$

We finish the section by proving the assertion in
 Remark \ref{remark 11.10.06}. Clearly, 
$$
2\sum_{\lambda\in\Lambda}(
\bar{\delta}_{\lambda}\varphi)
L^0_{\lambda}T_{\lambda}\varphi=
I_1+I_2, 
$$
with 
$$
I_1:=2\sum_{\lambda\in\Lambda_1}(\bar{\delta}_{\lambda}\varphi)
L^0_{\lambda}T_{\lambda}\varphi, 
\quad
I_2:=2\sum_{\lambda\in\Lambda_2}(\bar{\delta}_{\lambda}\varphi)
L^0_{\lambda}\varphi.
$$
We take $\tau_{\lambda}\equiv1$ on $\Lambda_{1}$
and 
notice that due to the symmetry of 
$\Lambda_1$ and the symmetry of $Dq_{\lambda}$ in $\lambda$ 
$$ 
I_1=
2\sum_{\lambda\in\Lambda_1}(\delta_{\lambda}\varphi)
L^0_{\lambda}\varphi
+
2h\sum_{\lambda\in\Lambda_1}(\delta_{\lambda}\varphi)
L^0_{\lambda}\delta_{\lambda}\varphi
$$
$$
=\sum_{\lambda ,\mu\in\Lambda_1}
(\delta_{\lambda}\varphi)
(\delta_{\lambda}q_{\mu})\Delta_{\mu}\varphi
+2\sum_{\lambda ,\mu\in\Lambda_1}
(\delta_{\lambda}\varphi)
(\delta_{\lambda}p_{\mu})\delta_{\mu}\varphi
$$
$$
+2\sum_{\lambda ,\mu\in\Lambda_1}
(\delta_{\lambda}\varphi)
(\delta_{\lambda}\chi_{\mu})\delta_{\mu}\delta_{\lambda}\varphi
=:I_1^{(1)}+I_1^{(2)}+I_1^{(3)}, 
$$
$$
I_2=I_2^{(1)}+I_2^{(2)},
$$
where in the notation $\xi=D\varphi/|D\varphi|$,
$\psi_{(\xi)}=\xi_{i}D_{i}\psi$,
$$
I_2^{(1)}=\tau_{0}^{2}\sum_{j=1}^{d}\sum_{ \mu\in\Lambda_1}
(D_{j}\varphi)
(D_{j}q_{\mu})\Delta_{\mu}\varphi=\tau_{0}^{2}|D\varphi|
\sum_{\mu\in\Lambda_{1}}q_{\mu(\xi)}\Delta_{\mu}\varphi,
$$
$$
I_2^{(2)}=
2 \tau_{0}^{2}\sum_{j=1}^{d}
\sum_{ \mu\in\Lambda_1}
(D_{j}\varphi)
(D_{j}p_{\mu})\delta_{\mu}\varphi=
2 \tau_{0}^{2}|D\varphi|
\sum_{ \mu\in\Lambda_1}
 p_{\mu(\xi)} \delta_{\mu}\varphi
$$
By Young's inequality,  
taking into account that 
$\chi_{\lambda}\geq\kappa/2>0$ for sufficiently small 
$h$, and that $c\geq c_0>0$, we have  
$$
I^{(j)}_1\leq(1/3) \sum_{\lambda\in\Lambda_1}
\cQ(\delta_{\lambda} \varphi)
+N\cQ(\varphi)\quad \text{for $j=1,3$}, 
\quad I^{(2)}_{1}\leq N\cQ(\varphi),
$$
$$
I^{(1)}_2\leq (1/3)\sum_{ \mu\in\Lambda_1}
\cQ(\delta_{\mu}\varphi)
+ \tau_{0}^{2}Nc_0^{-1}c
\sum_{\lambda\in\Lambda_2}|\delta_{\lambda}\varphi|^2, 
$$ 
$$
I^{(2)}_2\leq \tau_{0}^{2}Nc_0^{-1}c
\sum_{\lambda\in\Lambda}|\delta_{\lambda}\varphi|^2,
$$
where 
$N$ is a constant depending only on 
$\kappa$,  the number 
of elements in $\Lambda_1$ 
and on the supremum norm of 
the gradients of 
$p_{\lambda}$ and $q_{\lambda}$ in $x$. 
Summing up these inequalities 
and taking $\tau_{0}>0$ sufficiently 
small  we get 
\eqref{3.24.1} with $K_1=2N$,
  unit operator $\cK$, and with 
$\delta$ as close to $1$ as we wish.

\mysection{Proof of Theorem \protect\ref{theorem 4.7.2}} 
                                  \label{section 07.9.25.1}

If $Q=\bR^{d}$,   Theorem \ref{theorem 07.9.12.1} (i)  shows 
that equation \eqref{equation} 
has a unique bounded 
continuous solution $u$ for which
the partial derivatives  in $x\in\bR^d$
up to order $m$ are bounded continuous functions of $(t,x)$. 
However, the bounds, which can be extracted
from the proof of Theorem \ref{theorem 07.9.12.1}
for these  derivatives depend on the parameters $h$ and $T$. 
Our aim is to show  the existence of bounds, independent 
of $h$ and $T$ if $m=1$ and in addition to  
Assumption \ref{assumption 16.12.07.06}, 
 Assumptions \ref{assumption 1.26.11.06} and 
\ref{assumption 11.22.11.06}  
also hold. We will obtain such estimates by 
making use of the following version of 
the maximum principle.
It is probably worth drawing the reader's
attention to the fact that the assumption
that $c$ has certain sign is not used in
Lemma \ref{lemma 3.6.1}.

\begin{lemma}
                                        \label{lemma 3.6.1}
Let Assumption
\ref{assumption 16.12.07.06} with $m=0$ be satisfied and let
$\chi_{\lambda}\geq0$ for all $\lambda\in\Lambda_{1}$.
Let $v$ be a bounded
 function on $Q_{T}$, such that $v(\cdot,x)$
is measurable  for any $x\in Q$ and
the partial derivative $D_{t}v:=\partial v(t,x)/\partial t$
exists in $Q^{o}_{T}$. 
Let   
$F(t)\geq0$ be an integrable function on $[0,T]$,
and let $C(t,x)\geq0$ be a bounded function. Assume that
  for all $(t,x)\in Q^{o}_{T}$ we have
\begin{equation}
                                               \label{07.9.23.1}
D_{t}v \leq  L v + C\bar{v}_{+}+F,
\end{equation}
where $\bar{v} (t)=\sup\{v (t,x):x\in Q\}$.

Then in $[0,T]$ we have
\begin{equation}
                                        \label{3.11.1}
\bar{v}(t)\leq  G(t)e^{\nu t} 
+\int_{0}^{t} F(s)e^{\nu(t-s)}\,ds,
\end{equation}
where
$$
\nu:=\sup_{Q^{o}_{T}}(C-c),\quad
G(t)=\sup_{(s,y)\in
 \delta' Q_{t}}e^{-\nu s}v_{+}(s,y),
$$
$$
\delta'Q_{t}=(\{ 0 \}
\times Q)\cup\delta_{x}Q_{t}.
$$
\end{lemma}

Proof. 
First assume that $C=c=F=G=0$. 
In that case introduce $\tilde{v}(t,x)=
v(t,x)e^{Nt}$, where
$$
N=N_{h}=\sup_{Q_{T}}
h^{ -1}
\sum_{\lambda\in\Lambda_{1}}
\chi_{\lambda}.
$$
Observe that $v(t,x)=\tilde{v}(t,x)e^{-Nt}$,
$$
e^{-Nt}(D_{t}\tilde{v}(t,x)-N\tilde{v}(t,x))
\leq e^{-Nt}L\tilde{v}(t,x),
$$
$$
 D_{t}\tilde{v}(t,x) 
\leq  (L +N)\tilde{v}(t,x) 
$$
in $Q^{0}_{T}$, and in $Q_{T}$
$$
\tilde{v}(t,x)\leq J\tilde{v}(t,x):=I_{Q^{0}_{T}}\int_{0}^{t}
(L +N)(I_{Q_{T}}\tilde{v})(s,x)\,ds.
$$
Obviously, due to the choice
 of $N$ and the assumption that $\chi_{\lambda}\geq0$, 
if $u_{1}\geq u_{2}$ on $Q_{T}$, then
$Ju_{1}\geq Ju_{2}$ on $Q_{T}$. It follows that,
$\tilde{v}\leq J^{k}\tilde{v}$ for any $k$. 

Observe that for any bounded function $f$ 
$$
\sup_{Q}|Jf(t,\cdot)|\leq N'\int_{0}^{t}
\sup_{Q}|f (s,\cdot)|\,ds,
$$
where $N'$ is independent of $f$ and $t$.
It follows easily that  for any bounded function $f$ we have
$J^{k}f\to0$ uniformly on $Q_{T}$. Hence
$\tilde{v}\leq0$ and $v\leq0$ in $Q_{T}$. 

Now we consider another particular case
in which $G=c=0$ and $C=\text{const}\geq0$. Then
 observe that
the function 
$$
w(t,x)= \int_{0}^{t} (C\bar{v}_{+}+F)(s)\,ds
$$
satisfies 
$$
D_{t}w=  Lw+C\bar{v}_{+}+F.
$$
Therefore $\hat{v}:=v-w$
satisfies 
 $
D_{t}\hat{v}
\leq L \hat{v}
$. In addition $\hat{v}\leq0$ on $\delta'Q_{T}$. By the above,
on $Q_{T}$ we have $\hat{v}\leq0$, that is
$$
v(t,x)\leq  \int_{0}^{t}
(C\bar{v}_{+}+F)(s) \,ds\quad
\bar{v}_{+}(t)\leq \int_{0}^{t}
(C\bar{v}_{+}+F)(s) \,ds,
$$
and \eqref{3.11.1} follows by Gronwall's inequality 
and the fact that
$\nu=C$.

Now we allow $c\ne0$ and variable $C$ 
 but still assume that $G=0$.
In that  case take a large constant $M$ so that
$M>c$ and $M+\nu>0$ and for $\hat{v}(t,x)=v(t,x)e^{Mt}$
write
$$
e^{-Mt}(D_{t}\hat{v}(t,x)-M\hat{v}(t,x))
\leq e^{-Mt}L \hat{v}(t,x)+C(t,x)\bar{v}_{+}(t)+F(t).
$$
Dropping obvious values of arguments and introducing
$$
\Bar{\Hat v}(t)=\sup_{x\in Q}\hat{v}(t,x)\quad(=\bar{v}e^{Mt}),
$$
we find
$$
D_{t}\hat{v}  
\leq  L^{0} \hat{v}+(M-c)\hat{v} +C\Bar{\Hat v}_{+}+e^{Mt}F 
\leq  L^{0} \hat{v}+(M-c)\Bar{\Hat v}_{+}
$$
$$
 +C\Bar{\Hat v}_{+}
+e^{Mt}F \leq  L^{0} \hat{v}+(M+\nu)\Bar{\Hat v}_{+}+e^{Mt}F.
$$
It follows by the above that
$$
\bar{v}(t)e^{Mt}\leq  \int_{0}^{t}e^{Ms}F(s)
e^{(M+\nu)(t-s)}\,ds,
$$
which is equivalent to
\eqref{3.11.1}. By the way, notice that so far we have not
used the fact that $C\geq0$.

Now comes the general case in which we set
$$
w(t,x)=v(t,x)e^{-\nu t}-\beta,\quad
\beta=\sup_{(s,y)\in\delta'Q_{T}}v_{+}(s,y)e^{-\nu s}.
$$
Simple manipulations show that \eqref{07.9.23.1}
becomes
$$
D_{t}w\leq L^{0}w-c(w+\beta)-\nu(w+\beta)
+Ce^{-\nu t}\bar{v}_{+}+Fe^{-\nu t},
$$
where 
$$
e^{-\nu t}\bar{v}_{+}\leq\bar{w}_{+}+\beta,
\quad\bar{w}_{+}(t):=\sup_{Q}w_{+}(t,x)
$$
and, owing to the assumption that $C\geq0$,
the definition of $\nu$, and the fact that $\beta\geq0$,
$$
Ce^{-\nu t}\bar{v}_{+}
\leq C\bar{w}_{+}+C\beta\leq C\bar{w}_{+}+(c+\nu)\beta.
$$
It follows that
$$
D_{t}w\leq L^{0}w -(c+\nu)w
+C\bar{w}_{+}+Fe^{-\nu t},
$$
Since $w\leq0$ on $\delta'Q_{T}$
and $\sup(C-(c+\nu))=0$, by the above cases we have for
$t\leq T$ that 
$$
\bar{w}(t)\leq\int_{0}^{t}F(s)e^{-\nu s}\,ds,
$$ 
\begin{equation}
                                               \label{07.9.23.2}
\bar{v}(t)\leq e^{\nu t}\sup_{(s,y)\in\delta'Q_{T}}
v_{+}(s,x)e^{-\nu s}+e^{\nu t}\int_{0}^{t}F(s)e^{-\nu s}\,ds.
\end{equation}
We can put here $t=T$ and then, by using certain freedom
in choosing the end of the time interval, we can, actually,
set $T=t$ in \eqref{07.9.23.2}. Then we arrive at
\eqref{3.11.1} for all $t\leq T$.    
The lemma is proved.

\begin{corollary}
                                      \label{corollary 07.9.23.1}
Under the conditions of Lemma \ref{lemma 3.6.1}
if $\nu<0$, then
$$
\bar{v}(t)\leq\sup_{(s,y)\in\delta'Q_{t}}v_{+}(s,y)
+|\nu|^{-1}\sup_{[0,t]}F.
$$

\end{corollary}

{\bf Proof of Theorem \ref{theorem 4.7.2}}.  
Introduce
$$
V_{0}=u^{2},\quad
V_{1 }= \sum_{\lambda
\in\Lambda }|\bar{\delta}_{\lambda}u|^{2},
\quad 
\bar{V}_{k }(t)=\sup_{x\in Q}V_{k }(t,x),
$$
 recall that $F_{1}$ is introduced
in the statement of Theorem \ref{theorem 4.7.2} and set
$$
G=\sup_{\delta'Q_{T}}(|u|+\tau_{0}  |Du|
+U).
$$

By Corollary \ref{corollary 07.9.23.1}
applied to $u$ or $-u$ from the 
assumption that $c\geq c_{0} $ we obtain that
$$
V_{0 }=|u|^{2}\leq (G+c_{0}^{-1}F_{1})^{2}  .
$$
Now we
use the formula
$$
\delta_{ \lambda}(\psi\varphi)
=(\delta_{ \lambda}\psi)T_{ \lambda}\varphi
+\psi\delta_{ \lambda}\varphi
=(\delta_{ \lambda}\psi) \varphi
+\psi\delta_{ \lambda}\varphi
+h(\delta_{ \lambda}\psi) \delta_{ \lambda}\varphi ,\quad
\lambda\in\Lambda _{1},
$$
 to get
$$
L^{0} (\varphi^{2})=2\varphi L^{0} \varphi
+\cQ(\varphi).
$$
In particular,
$$
L^{0} V_{1}=2 \sum_{\lambda\in\Lambda}
(\bar{\delta}_{ \lambda}u)L^{0} \bar{\delta}_{ \lambda}u
+ 
\sum_{ \lambda\in\Lambda} 
\cQ 
(\bar{\delta}_{ \lambda}u).
$$
We observe that 
$$
L^{0} \bar{\delta}_{ \lambda}u
=\bar{\delta}_{ \lambda}L^{0} u
-L ^0 _{ \lambda}T_{ \lambda}u
$$
 if $\lambda\in\Lambda $, and 
in $Q^{o}_{T}$ 
$$
\bar{\delta}_{ \lambda}L^{0} u=
D_{t}\bar{\delta}_{ \lambda}u
+c\bar{\delta}_{ \lambda}u 
+ (\bar{\delta}_{ \lambda}c)T_{ \lambda}u
-\bar{\delta}_{ \lambda}f.
$$
Then in $Q^{o}_{T}$ we find
$$
L^{0} V_{1}-D_{t}V_{1}
-2cV_{1}=\sum_{\lambda\in\Lambda}
\cQ 
(\bar{\delta}_{ \lambda}u)
$$
$$
+2\sum_{\lambda\in\Lambda}
(\bar{\delta}_{ \lambda}u)( (\bar{\delta}_{ \lambda}c)
T_{\lambda}u- \bar{\delta}_{ \lambda}f)
-2\sum_{\lambda\in\Lambda}
(\bar{\delta}_{ \lambda}u)L^{0}_{ \lambda}
T_{ \lambda}u\,.
$$
We  use Assumption \ref{assumption 11.22.11.06} to conclude
\begin{equation}
                                          \label{07.9.25.1}
L^{0} V_{1}-D_{t}V_{1}
-2cV_{1}\geq
-K_{1} 
\cQ(u)-2(1-\delta)c\bar{V}_{1}+I,
\end{equation}
where
$$
I:=2\sum_{\lambda\in\Lambda}
(\bar{\delta}_{ \lambda}u)( (\bar{\delta}_{ \lambda}c)
T_{ \lambda}u- \bar{\delta}_{ \lambda}f).
$$
Notice that by Young's inequality
$$
2\sum_{\lambda\in\Lambda}
|(\bar{\delta}_{ \lambda}u)
( (\bar{\delta}_{ \lambda}c)T_{ \lambda}u|
\leq\delta^{2}\bar{V}_{1}+
 \bar{V}_{0}\sum_{\lambda\in\Lambda}
(\bar{\delta}_{ \lambda}c)^{2},
$$
where, for each $\lambda\in\Lambda_{1}$,
$$
\delta_{\lambda}c(t,x)=h^{-1}\int_{0}^{h}\lambda_{i}D_{i}c(
t,x+\lambda\theta)\,d\theta,
$$
$$
|\bar{\delta}_{\lambda}c(t,x)|^{2}\leq|\tau_{\lambda}
\lambda|^{2}
h^{-1}\int_{0}^{h}
|Dc(t,x+\lambda\theta)|^{2}\,d\theta
\leq|\tau_{\lambda}\lambda|^{2}\sup_{H_{T}}|Dc|^{2},
$$
so that
$$
\sum_{\lambda\in\Lambda}(\bar{\delta}_{ \lambda}c)^{2}
\leq\sup_{H_{T}}|Dc|^{2}
(1+|\Lambda_{1}|^{2}), 
$$
$$
2\sum_{\lambda\in\Lambda}
|(\bar{\delta}_{ \lambda}u)( 
(\bar{\delta}_{ \lambda}c)T_{ \lambda}u|
\leq\delta^{2}\bar{V}_{1}+
N\bar{V}_{0}\leq\delta^{2}\bar{V}_{1}+
N(F^{2}_{1}+G^{2}) .
$$
Similarly,
$$
2\sum_{\lambda\in\Lambda}
|(\bar{\delta}_{ \lambda}u) \bar{\delta}_{ \lambda}f|
\leq\delta^{2}\bar{V}_{1}+NF_{1}^{2}.
$$
  
Hence \eqref{07.9.25.1} yields
\begin{equation}
                                    \label{3.19.2} 
L^{0} V_{1 }-D_{t}V_{1 }
-2cV_{1 }\geq-2( \delta^{2}+ c-\delta c) \bar{V}_{1 }
-K_{1}\cQ(u)-N (F_{1}^{2}+G^{2}).
\end{equation}

Next, 
$$
L ^{0}(u^{2})-2cu^{2}-D_{t}(u^{2})
=2u(L ^{0}u-cu - 
D_{t}u)+ \cQ(u)
$$
$$
=-2uf+ \cQ(u) \geq-N(F_{1}^{2}+G^{2})
+ \cQ(u) .
$$

This, \eqref{3.19.2}, and the fact that 
$\delta^{2}+c-\delta c\geq0$
show that for
$$
W:=V_{ 1}+K_{1}u^{2}
,\quad\bar{W}=\sup_{Q}W(\cdot,x),
$$
we have
$$
L_{h}^{0}W-2cW-D_{t}W\geq
- 2( \delta^{2}+ c-\delta c) \bar{W}
-N(F_{1}^{2}+G^{2}).
$$
Now we want to use  
Corollary \ref{corollary 07.9.23.1}. Set
$$
\nu:=\sup_{Q^{o}_{T}}[2( \delta^{2}+ c-\delta c)-2c]
=2\delta\sup_{Q^{o}_{T}}(\delta-c)\leq
2\delta(\delta-c_{0}).
$$
If 
\begin{equation}
                                      \label{07.9.8.1}
\nu\leq-\delta c_{0}
\end{equation}
then
by Corollary \ref{corollary 07.9.23.1} we get 
$W\leq N(F_{1}^{2}+G^{2})$, 
which obviously implies \eqref{4.8.03}.
Finally,
observe that, 
if Assumption \ref{assumption 11.22.11.06}
is satisfied with a $\delta=\delta_{0}>0$, 
then it is also satisfied
with any $\delta\in(0,\delta_{0}]$ 
and the same $K_{1},\cK$.
Hence by modifying $\delta$ if necessary, 
so that $\delta\leq
c_{0}/2$, we satisfy  
\eqref{07.9.8.1} thus proving the theorem.

The following remark will be used in a subsequent paper
when we will be estimating higher order derivatives of $u$.
\begin{remark}
                                \label{remark 07.11.25.6}

Suppose that, instead of 
Assumption \ref{assumption 11.22.11.06}, Assumption
\ref{assumption 07.9.25.1} 
(see Section \ref{section 4.9.1})
is satisfied.
Then a simple inspection of the above
proof shows that in place of \eqref{3.19.2}  we would have
$$
L^{0} V_{1 }-D_{t}V_{1 }
-2cV_{1 }\geq-2( \delta^{2}+ c-\delta c) \bar{V}_{1 }
$$
$$ 
+ \delta
\sum_{ \lambda\in\Lambda}\cQ(\bar{\delta}_{ \lambda}u)
-K_{1}\cQ(u)-N(F_{1}^{2}+G^{2}).
$$
\end{remark}

\mysection{Some issues related to the
convergence $L_{h}\to\cL$}
                              \label{section 07.9.25.6}

There is a natural question about the relation
of the finite-difference operators $L=L_{h}$
with partial-differential operators. We certainly
want to apply the results of the present article
to investigating approximate solutions
of elliptic and parabolic second-order equations.
Then, given an elliptic operator
$$
\cL=a_{ij}D_{i}D_{j}+b_{i}D_{i}
$$
with variable coefficients,
 a natural question arises as to whether it is possible
to construct operators $L_{h}$ such that 
 they converge to $\cL$ and our 
assumptions are satisfied.

This question has little to do with the dependence of 
$a_{ij}$ and $b_{i}$ on $t$ and we assume that
$$
a_{ij}=a_{ij}(x),\quad b_{i}=b_{i}(x),
$$
and $a$ and $b$ are bounded and continuous along 
with their first-order derivatives.

It is not hard to see that under the symmetry assumption  (S)
the operators $L_{h}$ approximate $\cL$ with
\begin{equation}
                                             \label{07.9.20.1}
a_{ij}(x)=(1/2)
\sum_{\lambda\in\Lambda_{1}}
q_{\lambda}(x)\lambda_{i}\lambda_{j},
\end{equation}
$$
 b_{i}=
\sum_{\lambda\in\Lambda_{1}}
p_{\lambda}(x)\lambda_{i},
$$ 
 in the sense that $L_{h}\varphi\to
\cL\varphi$ as $h\downarrow0$ for all smooth $\varphi$.

A few basic examples describing  conditions
on $q_{\lambda}$, which guarantee that
our assumptions are satisfied, are given
in   Remarks \ref{remark 11.10.06}, 
\ref{remark 07.9.18.7}, and \ref{11.17.10.07}.
 However,
how this can be transformed into some conditions
in terms of $a_{ij}$ is not clear right away.

By the way, it is shown in \cite{DK2} that if $\cL$ admits
finite-difference approximations constructed by contracting
a fixed mesh and the approximating
operators satisfy the maximum principle, then they always
have the form ~\eqref{07.9.25.03} with $\Lambda_{1}=
-\Lambda_{1}$ and $q_{\lambda}=q_{-\lambda}$.
This form is nonunique and the issue of choosing
appropriate  $q_{\lambda}$ and $p_{\lambda}$ arises.

It is proved in \cite{Kr08} that, if the matrix $a$
is uniformly nondegenerate, then there always exist
$\Lambda_{1}$ and $q_{\lambda}$ possessing property (S),
such that \eqref{07.9.20.1} holds,
 $q_{\lambda}$ are as smooth as $a$ is,
and $q_{\lambda}\geq\kappa>0$,   
where $\kappa$ is a constant.

It is also proved  in \cite{Kr08} that if 
all values of the matrix $a$ lie in a closed convex
polyhedron in the set of nonnegative matrices and
$a(x)$ has {\em two\/} bounded derivatives, then again
there   exist
$\Lambda_{1}$ and $q_{\lambda}$
 possessing property (S),
such that \eqref{07.9.20.1} holds, and 
$\sqrt{q_{\lambda}}$
are Lipschitz continuous.

In  these two cases 
the issue of satisfying our assumptions
reduces to   representing  $b(x)$ appropriately.

There is a way to do so, used quite often 
in probabilistic
literature, by adding, if necessary, the set
$\Gamma=\{\pm e^1,...,\pm e^{d}\}$ to $\Lambda_{1}$,
where $\{e^{i}\}$ 
is the standard basis in $\bR^{d}$, defining
$$
p_{\pm e^{i}}=(b_{i})_{\pm},
$$
($t_{\pm}=(1/2)(|t|\pm t)$), 
and defining $p_{\lambda}=0$
on the remaining part of $\Gamma\cup\Lambda_{1}$.
There  is a certain inconvenience
 in this approximation,
which we discuss in the following example 
along with a way
to avoid it by using different $p_{\lambda}$'s.

\begin{example}
For $d=1$ consider the operator
$$
Lu=bDu ,
$$
where  $b=b(x)$  is a smooth function 
bounded along with its derivatives.
  
 If $b$ changes sign,
then, since in our setting $\chi_{\lambda}=hp_{\lambda}$
 is required to be  $\geq0$,
 we have to take $\Lambda_{1}$ consisting
of at least two points $\{\lambda_{1}, 
\lambda_{3}\}$. The most natural choice is
$\Lambda_{1}=\{\pm1\}$ and
$$
L_{h}\varphi(x)=P_{h}\varphi(x)
=b_{+}(x)\delta_{h,1}\varphi(x)
+b_{-}(x)\delta_{h,-1}\varphi(x). 
$$
Observe that on smooth $\varphi$ we have   
$L_{h}\varphi\to bD\varphi$
as $h\downarrow0$.

Notice that if $b$ changes sign,  
$p_{\pm1}=b_{\pm}$ are Lipschitz
continuous but need not be continuously differentiable
unless we impose a severe restriction on the behavior of
$b$ near the points where it vanishes. 
Actually, in this article the assumption that $q_{\lambda}$
and $p_{\lambda}$ are smooth can be replaced with the assumption
that they are Lipschitz continuous
  and then require \eqref{3.24.1}
be satisfied for each $t$ 
almost everywhere with respect to $x$
rather than for all $x$.
However,
in such case it is unrealistic to assume 
in Theorem \ref{theorem 4.7.2} that
$u$ is continuously differentiable in $x$. Generally, $u$
will be only Lipschitz continuous in $x$ and estimate 
\eqref{4.8.03} will hold almost everywhere rather that
everywhere in $Q_{T}$. More serious trouble occurs when
we want to estimate higher order derivatives, which we
will be concerned with in a subsequent article. Then we need
$p_{\lambda}$ to have higher order derivatives and this 
excludes many interesting cases.

  On the other hand, the reader may like
to check that with the above $p$ condition
\eqref{3.24.1} is satisfied (a.e.) 
with any $\delta\in(0,1)$
as long as $b$ is a {\em decreasing\/} function, 
which
agrees well with  the limit case of 
differential equations.

One can construct a different approximation
of $bD\varphi$ for which $p_{\lambda}$ 
are as smooth as $b$.
Indeed, take a constant $\theta$ 
such that $|b|+1\leq\theta$
and set
$$
p_{1}=b+\theta,\quad p_{-1}=\theta.
$$
Then $p_{\pm1}\geq1$,
again $p_{1}\delta_{h,1}\varphi+p_{-1}
\delta_{h,-1}\varphi\to bD\varphi$ 
on smooth $\varphi$,
and $p_{\pm1}$ are as smooth as $b$.
This method is somewhat close to adding an artificial
diffusion. However, we add it 
only to the finite-difference
approximation and not to the operator $bD$.

Now consider the operator
$$
Lu=(1/2)aD^{2}+bDu ,
$$
where we suppose that $a(x)\geq0$ and
$r:=\sqrt{a}$ and $b$ are one time 
differentiable with derivatives {\em uniformly\/}
continuous on $\bR$. Again take $\Lambda_{1}=\{\pm1\}$ an     
 construct
 $p_{\lambda}$ as in Remark \ref{remark 07.9.18.9} 
and define
$q_{\lambda}=a$ and $r_{\lambda}=r=\sqrt{a}$.

It is shown in  
Remark \ref{remark 08.1.5.1}   that
for $h$ sufficiently small,    Assumption 
\ref{assumption 11.22.11.06}   holds (perhaps with 
different $\delta$ and $K_{1}$),
  if 
$$
14(r')^{2}+b'\leq(1-\delta)c+K_{1}a.
$$ 
 
This condition describes what we need
from $\cL$ in the one-dimensional case
and it looks quite satisfactory. On the other hand,
it is yet  
stronger than the common assumption
$$
|r'|^{2}+b'\leq(1-\delta)c+K_{1}a,
$$
which along with other standard assumptions
guarantee that  solutions of $D_{t}u=\cL u-cu+f$
admit estimates of the {\em first\/}
 derivatives in $x$ independent
of the time interval.
 \end{example}

\mysection{Discussion of
Assumption \protect\ref{assumption 11.22.11.06}}
                                  \label{section 4.9.1}

In a subsequent paper about higher 
order derivatives estimates
we will impose the 
following assumption, which trivially
implies Assumption \ref{assumption 11.22.11.06}:
\begin{assumption}            \label{assumption 07.9.25.1} 
 We have $m\geq1$ and there
exist  a  constant    
$\delta \in(0,1] $ 
 and an operator $\cK=\cK_{h}\in\frK$,
 such that   
$$
2m \sum_{\lambda\in\Lambda }
(\bar{\delta}_{ \lambda}\varphi)
L^{0}_{ \lambda}T_{ \lambda}\varphi\leq
 (1-\delta)\sum_{\lambda \in\Lambda}
\cQ(\bar{\delta}_{ \lambda}\varphi)
$$
\begin{equation}
                                           \label{07.9.25.7}
+K_{1}\cQ(\varphi)
+
 2(1-\delta)c \cK \big(\sum_{\lambda\in\Lambda }
|\bar{\delta}_{ \lambda}\varphi|^{2}\big)  
\end{equation}
on $H_T$
for all smooth functions $\varphi $.
\end{assumption}

In this section
we are going to discuss Assumptions
\ref{assumption 11.22.11.06} and \ref{assumption 07.9.25.1}.
 Here we suppose that only
Assumptions \ref{assumption 16.12.07.06} and
\ref{assumption 1.26.11.06} are satisfied.

\begin{remark}
                                  \label{remark 07.9.18.6}

Condition \eqref{07.9.25.7} involves a mixture
of finite differences and derivatives. 
It is reasonable
to split it into  two parts, the combination
of which turns out to imply  \eqref{07.9.25.7}:
  For all smooth $\varphi$ we have on $H_{T}$
that
$$
2m\sum_{\lambda\in\Lambda_{1}}
(\bar{\delta}_{ \lambda}\varphi)
L^{0}_{ \lambda}T_{ \lambda}\varphi\leq
(1- \delta)\sum_{\lambda \in\Lambda_{1} }
\cQ(\bar{\delta}_{ \lambda}\varphi)
$$
\begin{equation}
                                           \label{4.9.3}
+K_{1}\cQ(\varphi)
+
2(1-\delta)c\cK \big(\sum_{\lambda\in\Lambda_{1}}
|\bar{\delta}_{ \lambda}\varphi|^{2}\big) , 
\end{equation}
that is \eqref{07.9.25.7}
holds with $\tau_{0}=0$ and
$$
2m \sum_{i=1}^{d}(D_{i}\varphi)
L^{0}_{ i} \varphi\leq
K_{1} \sum_{\lambda \in\Lambda_{1}}
\cQ(\bar{\delta}_{\lambda}\varphi)
$$
\begin{equation}
                                           \label{4.9.4}
+ (1/2)\delta c|D\varphi|^{2}+
 K_{1}  \cK (\sum_{\lambda\in\Lambda_{1}} 
|\bar{\delta}_{ \lambda}\varphi|^{2}) ,
\end{equation}  
where and below by $\cK$ we denote
generic operators (perhaps, depending on $h$)
of class $\frK$ and
$$
L_{ i}^0\varphi=\sum_{\lambda\in\Lambda_{1}} 
(h^{-1} D_{i}q_{\lambda} + D_{i}p_{\lambda})
\delta_{\lambda}\varphi.
$$

To show that \eqref{4.9.3} combined with
\eqref{4.9.4} imply \eqref{07.9.25.7}
if we choose 
small $\tau_{0}>0$ appropriately, 
observe that in terms of
$\bar{\delta}_{ \ell_{i}}=
\tau_{0}D_{i}$ equation 
\eqref{4.9.4} means that
$$
2m\sum_{\lambda\in\Lambda_{2}}
(\bar{\delta}_{\lambda}\varphi)
L^{0}_{\lambda} \varphi\leq
K_{1}\tau_{0}^{2}\sum_{\lambda\in\Lambda_{1}}\cQ
(\bar{\delta}_{\lambda}\varphi)
$$
\begin{equation}
                                           \label{4.9.5}
+(1/2)\delta c\sum_{\lambda\in\Lambda_{2}}
|\bar{\delta}_{ \lambda}\varphi|^{2}+
K_{1}\tau_{0}^{2}\cK
\big(\sum_{\lambda\in\Lambda_{1}} 
|\bar{\delta}_{ \lambda}\varphi|^{2}\big).
\end{equation}
By choosing $\tau_{0}$ 
so that $K_{1}\tau_{0}^{2}\leq
\delta/2$
and $K_{1}\tau_{0}^{2}
\leq  (1/2)\delta c_{0}$,
slightly redefining $\cK $ to absorb the 
second term on the right in \eqref{4.9.5}, and summing up
\eqref{4.9.5} and \eqref{4.9.3} we come
to an inequality which is
even somewhat stronger  than
\eqref{07.9.25.7} if $ \delta$ there is replaced
with $\delta/2$, which is irrelevant. 
\end{remark}
 
\begin{remark}
                                 \label{remark 07.10.29.1}

Assume that
 the symmetry condition (S) holds, $q_{\lambda}\geq0$,
$p_{\lambda}\geq0$, and $r_{\lambda}:=\sqrt{q_{\lambda}}$
are Lipschitz continuous in $x$ with a constant
independent of $t$. Then it turns out that
condition \eqref{4.9.4} is satisfied
with any $\delta\in(0,1)$,
$\tau_{\lambda}\equiv1$,
 and appropriate $K_{1}$ and unit $\cK$.

To show this   
observe that, for any unit $\xi$,
 $| q_{\lambda(\xi)}|\leq Nr_{\lambda}$
with $N$ being the doubled Lipschitz constant of $r_{\lambda}$.
In particular, by H\"older's inequality
\begin{equation}
                                         \label{07.11.26.1}
 \big(\sum_{\mu\in\Lambda_{1}}
q_{\mu(\xi)}\Delta_{\mu}\varphi\big)^{2}
\leq N
\big(\sum_{\mu\in\Lambda_{1}}
\sqrt{q}_{\mu }|\Delta_{\mu}\varphi|\big)^{2}
\leq N\sum_{\mu\in\Lambda_{1}}\cQ(\delta_{\mu}\varphi),
\end{equation}
which allows us to 
  make obvious changes in the estimates
of $I^{(j)}_{2}$ in the end of 
Section \ref{section 13.19.10.07},
one of the changes
 being that now we can allow $\cQ(\bar{\delta}_{\lambda}
\varphi)$ to enter the estimates with as large constant
as we wish. 
\end{remark}

In the following remark we discuss an estimate
which was crucial in the nonlinear setting
for establishing a rate of convergence of
difference approximations to the true solutions
of Bellman's equations in cylindrical domains
(see \cite{DK1}). We will see how using
different $\tau_{\lambda}$ can help.

\begin{remark}
                                    \label{remark 07.11.25.1}

Consider the situation when the coefficients
$q_{\lambda}$,  $p_{\lambda}$ and $c$, the free term
$f$, and the terminal data
$g$ also depend
on a parameter $y\in\bR $:
$$
q_{\lambda}=q_{\lambda}(t,z),\quad p_{\lambda}=p_{\lambda}
(t,z),\quad c=c(t,z),\quad f=f(t,z),\quad g=g(z) ,
$$
where $z=(x,y)\in\bR^{d+1}$.
Assume that these functions and their first derivatives
in $z$ are bounded on 
$H'_{T}=[0,T]\times\bR^{d+1}$
and continuous in $z$. Assume that 
 $c\geq c_{0}$ for all values of the arguments.
Suppose that \eqref{4.9.3} holds on $H_{T}$
 for $m=1$, $\tau_{\lambda}\equiv1$, any smooth $\varphi(x)$,
 and any value of the parameter
$y$ with $\cK$ perhaps depending on $y$ (as well as $h$
and $t$).
Assume also that
 the symmetry condition (S) holds, $q_{\lambda}\geq0$,
$p_{\lambda}\geq0$, and $r_{\lambda} =\sqrt{q_{\lambda}}$
are Lipschitz continuous in $z$ with a constant
independent of $t$.

Finally, suppose that in $H'_{T}$ we are given a bounded
function $u(t,x)=u(t,x,y)$
which satisfies \eqref{equation}
in $Q_{T}$ for each value of $y$. Of course, now in  
 \eqref{equation} we write $z$ in place of $x$. 

Take an $\varepsilon>0$ and set
$$
T_{\varepsilon}^{y}u (t,x,y)=u(t,x,y+\varepsilon),\quad
\delta_{\varepsilon}^{y} =\varepsilon^{-1}(
T_{\varepsilon}^{y} -1).
$$

We claim that in $[0,T]\times Q\times\bR$ it 
holds that
\begin{equation}
                                      \label{07.11.25.6} 
|\delta_{\varepsilon}^{y}u|
\leq N\big(\sup_{H_{T}\times\bR}
(|f|+|D_{z}f|)
+\sup_{(\delta'Q_{T})\times\bR }
(|u|+|\delta_{\varepsilon}^{y}u|
+\sum_{\lambda\in\Lambda_{1}}
|\delta_{\lambda}u|)\big),
\end{equation}
if $\varepsilon\in(0, \tau h]$, 
where $N$ and $\tau\in(0,1]$  depend 
 only on $\delta$, $c_{0}$, $K_{1}$,
the number of elements in $\Lambda_1$, $|\Lambda_{1}|$,
and the Lipschitz constants
of $c$, $r_{\lambda}$, and $p_{\lambda}$ with respect to $z
=(x,y)$.

To prove the claim observe that,
 although equation \eqref{equation} can be considered
in $Q_{T}$ as an equation with parameter $y$, we will
treat it  as an equation in $Q'_{T}=Q_{T}\times\bR$.
Then we denote by $\lambda_{0}$ the positive   vector
on the $y$-axis having the length $\varepsilon/h$,
take a $\tau>0$ to be specified later and introduce
$$
\Lambda_{1}'=\Lambda_{1}\cup\{\lambda_{0}\},\quad
q_{\lambda_{0}}=p_{\lambda_{0}}=0,\quad
\tau_{\lambda_{0}}=\tau h/\varepsilon,\quad
\tau_{\lambda}=1,\quad \lambda\in\Lambda_{1}.
$$

We now check that Assumption \ref{assumption 11.22.11.06} 
is satisfied for the new objects
 with $m=1$, $\tau_{0}=0$ and $\Lambda'_{1}$
in place of $\Lambda$. Owing to the assumption
that \eqref{4.9.3} holds, we immediately see
that the left-hand side of \eqref{3.24.1}
for the new objects is less than
$$
(1-  
\delta )\sum_{\lambda \in\Lambda_{1} }
\cQ( \delta _{ \lambda}\varphi)
+K_{1}\cQ(\varphi)
+
2(1-\delta)c\cK \big(\sum_{\lambda\in\Lambda_{1}}
| \delta _{ \lambda}\varphi|^{2}\big) +I , 
$$
where $I=I_{1}+I_{2}$,
$$
I_{1}= (\bar{\delta}_{\lambda_{0}}\varphi) 
\sum_{\mu\in\Lambda_{1}}(\bar{\delta}_{\lambda_{0}}q_{\mu})
\Delta_{\mu}T_{\lambda_{0}}\varphi,\quad
I_{2}=
2 (\bar{\delta}_{\lambda_{0}}\varphi) 
\sum_{\mu\in\Lambda_{1}}(\bar{\delta}_{\lambda_{0}}p_{\mu})
\delta_{\mu}T_{\lambda_{0}}\varphi.
$$
Since for smooth $\psi$,
$$
\bar{\delta}_{\lambda_{0}}\psi=
\tau\delta^{y}_{\varepsilon}\psi,\quad
|\bar{\delta}_{\lambda_{0}}\psi|\leq
\tau\sup|\partial\psi/\partial y|,
$$
we have that
$$
|I_{2}|\leq N\tau \cK\big(\sum_{\mu\in\Lambda_{1}'}|
\bar{\delta}_{\mu}\varphi|^{2}\big).
$$

Upon observing the following general properties
of finite-differences:
\begin{equation}
                                      \label{07.11.25.4}
 h\delta_{\lambda}\delta_{\mu}=
(T_{\mu}-1)\delta_{\lambda},\quad
\Delta_{\mu}T_{\lambda}=-\delta_{\mu}\delta_{-\mu}
+\delta_{\lambda}\delta_{\mu}
+\delta_{\lambda}
\delta_{-\mu}
\end{equation}
and combining them with
the estimate 
$$
\tau^{-1}|\bar{\delta}_{\lambda_{0}}q_{\mu}|=
|\delta_{\varepsilon}^{y}q_{\mu}|
=|2r_{\mu}\delta_{\varepsilon}^{y}r_{\mu}
+\varepsilon(\delta_{\varepsilon}^{y}r_{\mu})^{2}|
\leq N (\sqrt{q}_{\mu}+h), 
$$
and \eqref{07.11.26.1},
one easily shows that 
$$
|I_{1}|\leq N\tau\big(\sum_{\lambda\in\Lambda'_{1}}
\cQ(\delta_{\lambda}\varphi)+
\cK\big(\sum_{\mu\in\Lambda_{1}'}|
\bar{\delta}_{\mu}\varphi|^{2}\big)\big).
$$
Since $\varepsilon\leq\tau h$, we have that
$|\delta_{\lambda}\varphi|\leq|\bar{\delta}_{\lambda}\varphi|$
on $\Lambda_{1}'$ and  the above estimates show how to choose
$\tau\in(0,1]$ in order for 
Assumption \ref{assumption 11.22.11.06} to
be satisfied indeed.
 
By applying Theorem \ref{theorem 4.7.2}
 (and Remark \ref{remark 07.11.26.5})
we finish proving our claim.
The only point which is perhaps worth noting is that
by Theorem \ref{theorem 4.7.2} the constant $N$
from \eqref{07.11.25.6} also depends on $|\Lambda_{1}'|$.
However, $|\Lambda_{1}'|^{2}=|\Lambda_{1}|^{2} +\tau^{2}$.
\end{remark}

\mysection{Discussion of Assumptions
\protect\ref{assumption 11.22.11.06} and 
\protect\ref{assumption 07.9.25.1} in case that
$\tau_{\lambda}=1,\lambda\in\Lambda_{1}$}
                                  \label{section 07.11.26.1}

Here we suppose that only
Assumptions \ref{assumption 16.12.07.06} and  
\ref{assumption 1.26.11.06} are satisfied. Everywhere
below we set $\tau_{\lambda}=1$ 
for all $\lambda\in\Lambda_{1}$.

\begin{remark}  
                                     \label{remark 07.9.17.1}
Suppose that   the symmetry assumption
(S) is satisfied. 
Then the operators $L_{h}$ can be regarded as finite-difference
approximations of $\cL$ 
(see \eqref{07.10.11.1})  
in the sense that for any smooth $\varphi$
we have $L_{h}\varphi\to\cL\varphi$ as $h\downarrow0$.
If we are only interested in this property, then we can always assume 
that $p_{\lambda}\geq K_{1}$.

Indeed, if we do not have this inequality, then we 
take a sufficiently large
constant $K_{2}$ (independent of $h$),
redefine $p_{\lambda}$  
as $p_{\lambda}+K_{2}$.
 This will not violate the convergence
$L_{h}\varphi\to\cL\varphi$ since
$$
2\sum_{\lambda\in\Lambda_{1}}\delta_{\lambda}\varphi=
\sum_{\lambda\in\Lambda_{1}}[\delta_{\lambda}+\delta_{-\lambda}]\varphi
\to0
$$
if $\varphi$ is smooth.
\end{remark}

The following lemma is often used below and in the 
continuation of the present paper. 
\begin{lemma}                        \label{lemma 14.22.10.07}
Let $\alpha_{\lambda\mu}$ be a nonnegative function 
on $\bR^d$ 
for each 
$\mu\in\Lambda_1$ and $\lambda$ from a finite set 
of indices $\Lambda^{\prime}$. Assume that 
$$
\sum_{\mu\in\Lambda_1}
\sup_{\lambda\in\Lambda^{\prime}}\alpha_{\lambda\mu}\leq C
$$
for some function $C$ on $\bR^d$.
Then  there is a
$\cK\in\frK$ such that 
\begin{equation}                           \label{15.22.10.07}
h^2\sum_{\mu\in\Lambda_1,\lambda\in\Lambda^{\prime}}
\alpha_{\lambda\mu}(\delta_{\mu}f_{\lambda})^2
\leq
4C\cK
(\sum_{\lambda\in\Lambda^{\prime}}f^2_{\lambda}).
\end{equation}
on $\bR^{d}$ for any bounded Borel
function
$f=f_{\lambda}$ given on $\bR^d$.
\end{lemma}
Proof. Using 
$$
h^2(\delta_{\mu}f_{\lambda})^2=
((T_{\mu}-1)f_{\lambda})^2\leq 2(T_{\mu}f_{\lambda})^2
+2f_{\lambda}^2=2(T_{\mu}+1)f^2_{\lambda},
$$ 
we can estimate from above the left-hand side of 
\eqref{15.22.10.07}  by 
$$
2\sum_{\mu\in\Lambda_1}C_{\mu}
(T_{\mu}+1)\sum_{\lambda\in\Lambda^{\prime}}f_{\lambda}^2, 
$$
where 
$
C_{\mu}
:=\sup_{\lambda\in \Lambda^{\prime}}
\alpha_{\lambda\mu}. 
$ 
Hence we get \eqref{15.22.10.07} with 
$\cK\in\frK$ defined by 
$$ 
\cK(f)=\tfrac{1}{2C}
\sum_{\mu\in\Lambda_1}C_{\mu}(T_{\mu}+1)f.
$$
The lemma is proved. 

\begin{remark}
                                      \label{remark 07.9.18.5}

One can give sufficient conditions for \eqref{4.9.3}
to hold without involving test functions $\varphi$,
which makes them  more ``explicit"
and in combination
 with Remark \ref{remark 07.10.29.1}
covers many situations when
Assumption \ref{assumption 07.9.25.1}  is relatively
easy to check.
One set of these ``explicit" conditions
is given in this remark.
Here we also show why the operators of class $\frK$  
are useful and how the presence of $hp_{\lambda}$ 
in $\chi_{\lambda}$ entering the operator $\cQ$ 
on the right of \eqref{4.9.3} may help.  

Suppose that $\Lambda_{1}=-\Lambda_{1}$ and $q_{\lambda}=
q_{-\lambda}\geq0$  and set 
$r_{\lambda}=\sqrt{{q}_{\lambda}}$. 
Take a   $\delta\in(0,1/4)$ and
assume that on $H_{T}$
there are functions $r_{\lambda\mu}$, $p_{\lambda\mu}\geq0$,
$\lambda,\mu\in\Lambda_{1}$
such that 
\begin{equation}
                                            \label{07.9.17.2}
h^{2}(\delta_{\lambda}r_{\mu})^{2}
\leq   \delta (\chi_{\mu}+\chi_{\lambda})
+h^{2}r^{2}_{\lambda\mu},\quad\sum_{\mu\in\Lambda_{1}}
\sup_{\lambda\in\Lambda_{1}}r^{2}_{\lambda\mu}\leq
2\delta c,
\end{equation}
\begin{equation}
                                           \label{07.9.17.3}
h^{2}|\delta_{\lambda}p_{\mu}|\leq
\delta^{2}
 ( \chi_{\mu}+\chi_{\lambda}) 
+\delta h^{2}p_{\lambda\mu},
\quad\sum_{\mu\in\Lambda_{1}}
\sup_{\lambda\in\Lambda_{1}}p_{\lambda\mu}
\leq \delta c, 
\end{equation}
 By 
virtue of Remark \ref{remark 07.9.17.1} if $L_{h}$
are used for approximating $\cL$, 
we can change these operators
and have \eqref{07.9.17.2} and \eqref{07.9.17.3} satisfied
with $r_{\lambda\mu}=p_{\lambda\mu}=0$ 
for sufficiently small $h$,
provided that the Lipschitz
constants in $x$ of $r_{\lambda}$ and $p_{\lambda}$
are bounded with respect to $t$.

For a function $\xi_{\lambda}$ given on $\Lambda_{1}$ 
let us write
$$
|\xi|^{2}=\sum_{\lambda\in\Lambda_{1}}|\xi_{\lambda}|^{2}.
$$
Then it turns out that condition \eqref{4.9.3}  
is satisfied 
if on $H_T$ for  all 
functions $\xi_{\lambda}$  
we have
$$
10m^{2}(1-4\delta)^{-1}
J_{ 1}
+2m^{2}(1-4\delta)^{-1}J_2
$$
$$
+2\delta m^{2}
  \sum_{\lambda,\mu \in\Lambda_{1}} \xi_{\lambda} ^{2}
|\delta_{\lambda}p_{\mu}|  
+2m\sum_{ \lambda,\mu\in\Lambda_{1}}
\xi_{\lambda}\xi_{\mu}(\delta
_{\lambda}p_{\mu}+(\delta_{\lambda}r_{\mu})^{2})
$$
\begin{equation}
                                             \label{07.9.16.1}
\leq(2 -8\delta)c|\xi|^{2}
+K_{1}\sum_{ \lambda \in\Lambda_{1}}
 \xi_{\lambda} ^{2}\chi_{\lambda}
+ \delta h^{-2}\sum_{\lambda\in\Lambda_{1}}
\chi_{\lambda}|\xi_{\lambda}+\xi_{-\lambda}|^{2},
\end{equation}
where
$$
J_{ 1}=\sum_{\mu,\lambda\in\Lambda_{1}}
\xi_{\lambda}^{2}(\delta_{\lambda}r_{\mu})^{2}, 
\quad
J_2=
\sum_{\mu\in\Lambda_{1}}\big(
\sum_{\lambda\in\Lambda_{1}}\xi_{\lambda}\delta_{\lambda}
r_{\mu}\big)^{2}.
$$
 
To prove this, use  formulas \eqref{07.11.25.4}.
Also drop the summation sign having repeated indices
in $\Lambda_{1}$ 
to see that
$$
2
(\delta_{\lambda}\varphi)Q_{\lambda}T_{\lambda}\varphi
=
(\delta_{\lambda}\varphi)
(\delta_{\lambda} q_{\mu})(\delta_{\lambda}
\delta_{\mu}+\delta_{\lambda}\delta_{-\mu}
+ \Delta_{\mu})\varphi
$$
$$
 =2\xi_{\lambda} (\delta_{\lambda} q_{\mu})\delta_{\lambda}
\delta_{\mu}\varphi
-
\xi_{\lambda}(\delta_{\lambda} q_{\mu})
\delta_{\mu}\delta_{-\mu} \varphi
= I_{11}+ I_{12}+ I_{21}+I_{22},
$$
where   $\xi_{\lambda}=\delta_{\lambda}\varphi$ and
$$
I_{11}= 4\xi_{\lambda}
(\delta_{\lambda}r_{\mu})r_{\mu}\delta_{\lambda}\delta_{\mu}
\varphi,
\quad
 I_{12}= 2
h\xi_{\lambda}(\delta_{\lambda}r_{\mu})^{2}
 \delta_{\lambda}\delta_{\mu}\varphi,
$$
$$
I_{21}=-2
\xi_{\lambda}
(\delta_{\lambda}r_{\mu})r_{\mu}\delta_{-\mu}\delta_{\mu}
\varphi,
\quad
I_{22}=  h\xi_{\lambda}(\delta_{\lambda}r_{\mu})^{2}
\Delta_{\mu}\varphi=
2\xi_{\lambda}(\delta_{\lambda}r_{\mu})^{2}\xi_{\mu}. 
$$
Before starting to estimate $I_{ij}$, we note that
as $h\downarrow0$ the terms $I_{11}$, 
$I_{22}$, and $I_{12}$
disappear if $\varphi$ is twice continuously
differentiable due to the symmetry of $\Lambda_{1}$.
 In that case there is no
need to estimate them.
For fixed $h$ they are present and estimating them
is only possible under stronger assumptions than
in the case of partial differential equations.

Now notice that by Young's inequality
$$
mI_{11}\leq(1/2)(1-4\delta)I
+ 8m^{2}(1-4\delta)^{-1}J_{ 1},
$$
where
$$
I= \sum_{\lambda\in\Lambda_{1}}\cQ(
\delta_{\lambda}\varphi)
=\sum_{\lambda,\mu\in\Lambda_{1}}\chi_{\mu}
| \delta_{\mu}\delta_{\lambda}\varphi|^{2}.
$$
Similarly,  
$$
mI_{21}\leq(1/2)(1-4\delta)I+
2m^{2}(1-4\delta)^{-1}
 J_2.
$$

Next, owing to \eqref{07.9.17.2}
$$
mI_{12}=2m \big(\xi_{\lambda}\delta_{\lambda}r_{\mu}\big)
\big(h(\delta_{\lambda}r_{\mu})\delta_{\lambda}
\delta_{\mu}\varphi\big)
$$
$$
\leq 2m^{2}J_{1}
+(1/2)
h^{2}\sum_{\lambda,\mu\in\Lambda_{1}}
(\delta_{\lambda}r_{\mu})^{2}(\delta_{\lambda}\delta_{\mu}\varphi)^{2}
$$
$$
\leq 2m^{2}J_{1}+\delta I
+(1/2)h^{2}
\sum_{\lambda,\mu\in\Lambda_{1}}r_{\lambda\mu}^{2}
(\delta_{\lambda}\delta_{\mu}\varphi)^{2},
$$
where the last term 
by virtue of 
Lemma \ref{lemma 14.22.10.07} is  
estimated by
$$
4\delta c\cK(\sum_{\lambda \in\Lambda_{1}}
( \delta_{\lambda}\varphi)^{2}). 
$$
By collecting the above estimates we obtain
$$
2m(\delta_{\lambda}\varphi)Q_{\lambda}T_{\lambda}\varphi
\leq (1-3\delta)I
+
10m^{2}(1-4\delta)^{-1}J_{1}
$$
\begin{equation}
                                           \label{07.9.18.1}
+2m^{2}(1-4\delta)^{-1}J_{2}
+2m(\delta_{\lambda}r_{\mu})^{2}
\xi_{\lambda}\xi_{\mu}
+4\delta c\cK
(\sum_{\lambda \in\Lambda_{1}}
( \delta_{\lambda}\varphi)^{2}).
\end{equation}
Next,
$$
2(\delta_{\lambda}\varphi)P_{\lambda}T_{\lambda}\varphi
=2\xi_{\lambda}(\delta_{\lambda}p_{\mu})\xi_{\mu}+
2h\xi_{\lambda}(\delta_{\lambda}p_{\mu})\delta_{\lambda}
\delta_{\mu}\varphi,
$$
where the last term is majorated by
$$
2\delta
m\sum_{\lambda,\mu\in\Lambda_{1}}\xi_{\lambda}^{2}
|\delta_{\lambda}p_{\mu}| +
(1/2)
m^{-1}\delta^{-1}h^{2}\sum_{\lambda,\mu\in\Lambda_{1}} 
|\delta_{\lambda}p_{\mu}|(
\delta_{\lambda}
\delta_{\mu}\varphi)^{2}.
$$
We use assumption \eqref{07.9.17.3} and proceed as
while estimating $I_{12}$. Then we see that
$$
2m(\delta_{\lambda}\varphi)P_{\lambda}T_{\lambda}\varphi
\leq  2m\xi_{\lambda}(\delta_{\lambda}p_{\mu})\xi_{\mu}
$$
\begin{equation}
                                           \label{07.9.18.2}
+2\delta
m^{2}\sum_{\lambda,\mu\in\Lambda_{1}}\xi_{\lambda}^{2}
|\delta_{\lambda}p_{\mu}|+\delta I +2\delta
c\cK(\sum_{\lambda \in\Lambda_{1}} (
\delta_{\lambda}\varphi)^{2}).
\end{equation}
Finally, upon combining \eqref{07.9.18.2}
with \eqref{07.9.18.1} we obtain
$$
2m(\delta_{\lambda}\varphi)L^{0}_{\lambda}T_{\lambda}\varphi
\leq (1-2\delta)I+
10m^{2}(1-4\delta)^{-1}
J_{1}
 +2m^2(1-4\delta)^{-1}J_2 
$$
$$
+2\delta
m^{2}\sum_{\lambda,\mu\in\Lambda_{1}}\xi_{\lambda}^{2}
|\delta_{\lambda}p_{\mu}|
+6\delta
c\cK(\sum_{\lambda \in\Lambda_{1}} (
\delta_{\lambda}\varphi)^{2})
$$
\begin{equation}                          \label{07.11.26.8}
+2m\xi_{\lambda}
 ((\delta_{\lambda}r_{\mu})^2+
(\delta_{\lambda}p_{\mu}))\xi_{\mu}.
\end{equation}
We use the fact that
$h^{-1}(\delta_{\lambda}+\delta_{-\lambda})=-
\delta_{\lambda}\delta_{-\lambda}$, use
 assumption \eqref{07.9.16.1},
and  take into account that for any 
$\cK^{\prime}\in\frK$ 
$$
(2-8\delta)
c\sum_{\lambda\in\Lambda_1}|\delta_{\lambda}\varphi|^2+
6\delta c\cK^{\prime}
(\sum_{\lambda\in\Lambda_1}|\delta_{\lambda}\varphi|^2)
=(2-2\delta)c
\cK(\sum_{\lambda\in\Lambda_1}|\delta_{\lambda}\varphi|^2)
$$
with an appropriate $\cK\in\frK$.
Then
 we estimate 
the right-hand side of \eqref{07.11.26.8} by 
$$
(1-\delta)I+K_{1}\cQ(\varphi)  
+(2-2\delta)c\cK \big(
\sum_{\lambda\in\Lambda_{1}}(\delta_{\lambda}\varphi)^{2}\big),  
$$
and we see that 
\eqref{4.9.3} is satisfied
indeed.
\end{remark}

\begin{remark}
It is easy to see that if 
\eqref{07.9.17.2} and \eqref{07.9.17.3} hold and we assume
that  inequality  \eqref{07.9.16.1} is satisfied with 
an additional term $\Psi(t,x,\xi)$ on its 
left-hand side for some function 
$\Psi$ of $t$, $x$ and
$\xi=(\xi_{\lambda})_{\lambda\in\Lambda_1}$, 
then inequality \eqref{4.9.3} holds with the additional
term 
$\Psi(t,x,(\delta_{\lambda}\varphi)_{\lambda\in\Lambda_1})$ 
on its left-hand side.
\end{remark}

\begin{remark}
                                      \label{remark 07.9.18.7}
Assume (S), assume that $q_{\lambda}\geq0$,
$p_{\lambda}\geq0$ 
 and let $m\geq2$. 
Then it turns out that  Assumption \ref{assumption 11.22.11.06} 
is satisfied for $\delta=1/10$  
and appropriate $\tau_{0}>0$ or $\tau_{0}=0$ 
if $c_{0}$ is sufficiently large
(independently of   $h$).

Indeed, it is well known that the Lipschitz constant
of the square root of a nonnegative twice
continuously
differentiable function $w(x)$ is controlled by 
the supremums of 
its second order derivatives.
Therefore,  Remark \ref{remark 07.9.18.5} 
(where we take $r_{\lambda\mu}=\delta_{\lambda}r_{\mu}$
and $p_{\lambda\mu}=10|\delta_{\lambda}p_{\mu}|$)
immediately implies
that condition \eqref{4.9.3} is satisfied 
if $c_{0}$ is large enough.

 That condition \eqref{4.9.4} is satisfied
follows from Remark \ref{remark 07.10.29.1}
and again from the fact  that for twice continuously
differentiable $w$  
on $\bR^{d}$
we have
$$
|D w|^{2}\leq  2w\sup_{x\in\bR^{d},|\xi|=1}|w_{(\xi)(\xi)}(x)|.
$$

At this point we do not even need large $c_{0}$.
Referring to Remark \ref{remark 07.9.18.6}
we obtain what we have claimed.

Actually, above in this remark we used that $m\geq2$
only to guarantee that $r_{\lambda}$ are Lipschitz
continuous in $x$ with a constant $N'$ independent of $t$.
If we just assumed this last property, then our argument
about \eqref{4.9.3} would become even shorter. In addition,
what was said about \eqref{4.9.4} is still valid.

\end{remark}

\begin{remark}  
                                     \label{remark 3.8.1}
 
Under the symmetry assumption (S)
and the assumption that $q_{\lambda}\geq0$
and $p_{\lambda}\geq0$
one can give a rougher condition without using $\xi_{\lambda}$
and implying \eqref{07.9.16.1}
with $m=1$ and sufficiently small $\delta$.
Then  \eqref{4.9.3}  
will be satisfied as long as conditions \eqref{07.9.17.2} and  
\eqref{07.9.17.3} are.

By the way, 
recall that, for all small $h$, 
one can always satisfy
  conditions \eqref{07.9.17.2} and   
\eqref{07.9.17.3}   on  the account of
modifying if necessary $p_{\lambda}$ 
if the Lipschitz
constants of $r_{\lambda}$ in $x$ are bounded in $t$
(see Remark \ref{remark 07.9.18.5}). 

 By the inequality  
$$
\sum_{ \lambda,\mu\in\Lambda_{1}}\eta_{\lambda}\eta_{\mu} 
\alpha_{\lambda\mu}\leq
\sum_{\lambda\in\Lambda_{1}}\eta_{\lambda}^{2}
\sum_{\mu\in\Lambda_{1}}|\alpha_{\lambda\mu}
+\alpha_{\mu\lambda}|
$$
  we have 
$$
\sum_{\lambda,\mu\in\Lambda_1}
\xi_{\lambda}\xi_{\mu}\delta_{\lambda}p_{\mu}
\leq \sum_{\lambda\in\Lambda_1}\xi_{\lambda}^2
\sum_{\mu\in\Lambda_1}
|\delta_{\lambda}p_{\mu}+\delta_{\mu}p_{\lambda}|, 
$$
$$
J_2=\sum_{\lambda,\nu,\in\Lambda_1}
\xi_{\lambda}\xi_{\nu}
\sum_{\mu\in\Lambda_1}
(\delta_{\lambda}r_{\mu})(\delta_{\nu} r_{\mu})
\leq 2\sum_{\lambda\in\Lambda_1}\xi_{\lambda}^2
\sum_{\nu\in\Lambda_1}|\sum_{\mu\in\Lambda_1}
(\delta_{\lambda}r_{\mu})\delta_{\nu} r_{\mu}|. 
$$
 
A simple argument based on 
 the above estimates 
and continuity shows that \eqref{07.9.16.1}
holds with $m=1$ and a small $\delta>0$ if
for any $\lambda\in\Lambda_{1}$
$$
10
\sum_{\mu\in\Lambda_{1}}(\delta_{\lambda}r_{\mu})^{2}+
4\sum_{\nu\in\Lambda_1}|\sum_{\mu\in\Lambda_1}
(\delta_{\lambda}r_{\mu})\delta_{\nu} r_{\mu}|
$$
\begin{equation}
                                    \label{07.10.29.01}
+2
\sum_{\mu\in\Lambda_{1}}|\delta_{\lambda}
p_{\mu}+\delta_{\mu}p_{\lambda}
+(\delta_{\lambda} r_{\mu})^{2}+(\delta_{\mu} r_{\lambda})^{2}|\leq
c+K_{1}q_{\lambda}.
\end{equation}

Condition \eqref{07.10.29.01} 
basically means that if for a $\lambda\in\Lambda_{1}$
at some point in $H_T$ the value $q_{\lambda}$
is small, then either $\delta_{\lambda}r_{\mu}$,
$\delta_{\mu} r_{\lambda}$,
$\delta_{\mu}p_{\lambda}$,
and $\delta_{\lambda}p_{\mu}$ should be small or $c$
be large at this point. As the point varies, the
dominating terms may change roles.

\end{remark}

\begin{remark}  
                                     \label{remark 07.9.18.9}
There are cases when it is preferable
to keep the last term on the left in \eqref{07.9.16.1}
as is. 

To see a reason for that, let $m=d=1$, take 
a constant $\theta\geq1$ and define
$\Lambda_{1}=\{\pm 1\}$, 
$$
q_{\lambda}\equiv0,\quad
p_{1}(t,x)=(1/2)b(x)+\theta,\quad 
p_{-1}=-(1/2)b(x)+\theta, 
$$ 
where $b(x)$ is a {\em decreasing\/} function
with bounded derivative
such that $|b|\leq1$. 
Observe that for any $\theta$ we have $L^{0}_{h}\varphi
\to b\varphi'$ as $h\downarrow0$ if $\varphi$ is smooth.

Now   notice that
condition \eqref{07.9.17.2} is trivially satisfied since $r_
{\lambda}\equiv0$. Condition \eqref{07.9.17.3}
is satisfied with any $\delta>0$ and
$p_{\lambda\mu}=0$ if $h\leq1$ and
$\theta\geq100\sup|b'|+1$ since 
$$
\chi_{\lambda}=hp_{\lambda}
\geq h(\theta-|b|)\geq h(\theta-1)
$$

The left-hand side of 
\eqref{07.9.16.1} is  one half of
$$
 4\delta \sum_{i=\pm1}\xi_{i}^{2} |\delta_{i}b|
+ (\xi_{1}-\xi_{-1})^{2}[\delta_{1}b-\delta_{-1}b]+
 (\xi_{1}^{2}-\xi^{2}_{-1}) (\delta_{1}b+\delta_{-1}b).
$$
Here the middle term is  
nonpositive since $b$ is decreasing.
Also the last term is majorated by
$$
N|\xi_{-1}+\xi_{1}|\,|\xi_{1}-\xi_{-1}|
\leq Nh(\xi_{1}^{2}+\xi_{-1}^{2})+Nh^{-1}(\xi_{-1}+\xi_{1})^{2}.
$$

Furthermore, concerning the right-hand side of \eqref{07.9.16.1}
observe that
$$
\sum_{\lambda}\xi^{2}_{\lambda}\chi_{\lambda}
\geq h(\theta-1)\sum_{\lambda}\xi^{2}_{\lambda},\quad
h^{-2}\sum_{\lambda}\chi_{\lambda}(\xi_{\lambda}+\xi_{-\lambda})^{2}
\geq h^{-1}(\theta-1)(\xi_{1}+\xi_{-1})^{2}.
$$
It follows easily that, no matter how small $c_{0}$
is, for sufficiently small $\delta$
and large $\theta$ condition \eqref{07.9.16.1}
and, by Remark \ref{remark 07.9.18.5}, condition
\eqref{4.9.3} are satisfied. This along with the almost 
obvious fact that \eqref{4.9.4} holds shows that
Assumption \ref{assumption 11.22.11.06} is satisfied
as well.

To finish the remark notice that if we tried to check condition
\eqref{07.10.29.01} we would fail to do that
for small $h$ no matter how large $\theta$ is
unless $c $ is large enough.

\end{remark}
         
\begin{remark}
                                       \label{remark 08.1.5.1} 
 We
continue the analysis of the one-dimensional situation started in
Remark
\ref{remark 07.9.18.9}. So, we assume that $d=m=1$ and we have in
mind approximating an operator $\cL\varphi(x)=(1/2)a(x)\varphi''(x)
+b(x)\varphi'(x)$. As in Remark \ref{remark 07.9.18.9}
we assume that $|b|\leq1$ and, in addition,
assume that $a\geq0$ and $r:=\sqrt{a}$ and $b$ are one time 
differentiable with derivatives {\em uniformly\/}
continuous on $\bR$.

As in Remark \ref{remark 07.9.18.9}
take $\Lambda_{1}=\{\pm1\}$ and define $p_{\pm1}$
for an appropriate $\theta$. Then also set
$r_{\mu}=r$. Now observe that both parts of 
\eqref{07.9.16.1} are order-two
homogeneous functions of $\xi_{\lambda}$.
Therefore, it suffices to check
\eqref{07.9.16.1} assuming that
\begin{equation}
                                            \label{08.1.5.4}
\sum_{\lambda\in\Lambda_{1}}\xi^{2}_{\lambda}=1.
\end{equation}

Then, due to the assumption that $r'$ and $b'$
are uniformly continuous on $\bR$, it is not hard to see that
$$
J_{1}\sim 2(r')^{2},\quad J_{2}\sim2(r')^{2}(\xi_{1}
-\xi_{-1})^{2},
$$
$$
\sum_{\lambda,\mu\in\Lambda_{1}}
\xi_{\lambda}\xi_{\mu}(\delta_{\lambda}r_{\mu})^{2}
\sim(r')^{2}(\xi_{1}+\xi_{-1})^{2},
$$
where by $\alpha\sim\beta$ we mean that for small $h$
the difference $|\alpha-\beta|$ can be absorbed
into $c\geq c_{0}>0$ with as small coefficient
as we wish. Therefore, upon recalling the estimates
from  Remark \ref{remark 07.9.18.9}, we see that condition
\eqref{07.9.16.1} is satisfied for all small $h$ if
$$
20(1-4\delta)^{-1}(r')^{2}+4(1-4\delta)^{-1} 2
(r')^{2}(\xi_{1}-\xi_{-1})^{2}
+2(r')^{2}(\xi_{1}+\xi_{-1})^{2}
$$
\begin{equation}
                                            \label{08.1.5.3}
+2\delta|b'|+(\xi_{1}-\xi_{-1})^{2}b'
\leq(2-9\delta)c+K_{1}a+\delta h^{-1}(\theta-1)(\xi_{1}+\xi_{-1}
)^{2}.
\end{equation}
On the account of (assumption \eqref{08.1.5.4} and)
 the presence of
$h^{-1}(\theta-1)(\xi_{1}+\xi_{-1})^{2}$ on the right 
of \eqref{08.1.5.3}, it suffices to check \eqref{08.1.5.3}
for small $h$ assuming that the inequality
$|\xi_{1}+\xi_{-1}|\leq h^{1/2}$ holds. It follows that
\eqref{08.1.5.3} holds for small $h$ if
it holds with $(2-10\delta)c$ in place of 
$(2-9\delta)c$ but only for $\xi_{\lambda}$
satisfying $\xi_{1}=-\xi_{-1}$. In that case
$(\xi_{1}-\xi_{-1})^{2}=4\xi_{1}^{2}=2$ and 
\eqref{08.1.5.3} holds if
$$
28(1-4\delta)^{-1}(r')^{2}+2\delta|b'|+2b'
\leq(2-10\delta)c+K_{1}a.
$$
Since we would be satisfied if \eqref{07.9.16.1}
held with at least one $\delta>0$, we see that,
under the assumption of the present remark,
\eqref{07.9.16.1} (perhaps with different $\delta$ and $K_{1}$) is
indeed satisfied for small $h$
 if 
$$
14(r')^{2}+b'\leq(1-\delta)c+K_{1}a.
$$
\end{remark}

\begin{remark}  
                                     \label{remark 07.9.19.1}
There are multi-dimensional analogs  
of the situation
in Remark \ref{remark 07.9.18.9}.
For instance, let $U(x)$ be a concave function
with bounded derivatives and assume that
$$
p_{\mu}=(DU,\mu)+\theta,
$$
  where 
$\theta\geq
 \max_{\mu\in\Lambda_1}|\mu| \sup|DU|$, 
so that $p_{\mu}\geq0$.
 Then for small $h$ we have
$$
\sum_{\lambda,\mu\in\Lambda_{1}}\xi_{\lambda}\xi_{\mu}\delta
_{\lambda}p_{\mu}\sim\sum_{\lambda,\mu\in\Lambda_{1}}
\xi_{\lambda}\xi_{\mu}(D^{2}U\mu,\lambda)
=(D^{2}U\eta,\eta)\leq0,
$$
where $\eta=\sum\xi_{\lambda}\lambda$.  

\end{remark}
\begin{remark}  
                                     \label{remark 3.8.01}
Recall that $\delta_{\lambda}=\delta_{h,\lambda}$
and assume that \eqref{4.9.4} holds for all
small $h>0$ with $\cK$ perhaps depending on $h$.
Then it turns out that 
for   all $t\in[0,\infty]$ 
\begin{equation}
                                         \label{3.9.1}
\sum_{\lambda\in\Lambda_{1}}\lambda q_{\lambda}(t,x )
\quad\hbox{is independent of}\quad  
 x .
\end{equation}

To show this observe that, since the values of the 
first derivatives of $\varphi$ at a fixed point
have nothing to do with the increments of $\varphi$,
\eqref{4.9.4} is equivalent to saying that
$$
\sum_{i=1}^{d}(L ^{0}_{ i}\varphi)^{2}
\leq N(\sum_{\lambda,\nu\in\Lambda_{1}}
\chi_{\lambda}|\delta_{ \lambda}\delta_{ \nu}\varphi|^{2}
+\cK (
 \sum_{\lambda \in\Lambda_{1}}
 |\delta_{ \lambda} \varphi|^{2})),
$$
where the constant $N$ can be easily computed
given $\sup c$ and $K_{1}$.
It follows that
\begin{equation}
                                        \label{07.8.12.3}
\sum_{i=1}^{d}(L ^{0}_{ i}\varphi)^{2}
\leq N(\sum_{\lambda,\nu\in\Lambda_{1}}
\chi_{\lambda}|\delta_{ \lambda}\delta_{ \nu}\varphi|^{2}
+\sup_{\bR^{d}}
 \sum_{\lambda \in\Lambda_{1}}
 |\delta_{ \lambda} \varphi|^{2}).
\end{equation}

Now  
multiply \eqref{07.8.12.3} by $h^{ 2}$ and let  $h
\downarrow0$. Then we obtain
$$
\sum_{i=1}^{d}\big[(D_{j}\varphi(x))\sum_{\mu\in\Lambda}
D_{i}
q_{\mu}(t,x)\mu^{j}\big]^{2}\leq N\sup_{\bR^{d}}|\varphi|.
$$
which leads to the conclusion that for smooth $\varphi$, 
$i=1,...,d$, and all
$\lambda\in\Lambda$
$$
(D_{j}\varphi(x))\sum_{\mu\in\Lambda}
D_{i}
q_{\mu}(t,x)\mu^{j}=0.
$$
This is equivalent to saying that \eqref{3.9.1}
holds.
\end{remark}

\begin{remark}
                                  \label{remark 3.26.1}
Additionally to 
Assumptions  \ref{assumption 16.12.07.06} and 
\ref{assumption 1.26.11.06} suppose that  
\begin{equation}
                                               \label{3.25.4}
 \lambda+\nu\not\in\Lambda_{1}\quad
\forall \lambda,\nu\in\Lambda_{1}.
\end{equation}
 It turns out that in this case
\eqref{4.9.4}
is satisfied for all small $h$ 
only if for any $\lambda\in\Lambda_{1}$
 
(i) either $-\lambda\in\Lambda_{1}$ and
$q_{\lambda}(t,x)=q_{-\lambda}(t,x)+r_{\lambda}(t)$
for a function $r_{\lambda}(t)$
  independent of $x$, 

(ii) or $-\lambda\not\in\Lambda_{1}$ and
 $q_{\lambda}$ is independent of $x$.

In particular,
\begin{equation}
                                              \label{3.26.2}
Q_{\nu}\varphi=(1/2)\sum_{\lambda\in\Lambda_{1}
\cap(-\Lambda_{1})}
(\delta_{\nu}q_{\lambda})\Delta_{\lambda}\varphi
\quad\quad(\sum_{\emptyset}...:=0).
\end{equation}
 
We may concentrate on proving our claim  
assuming that
$q_{\lambda}$ is independent of $t$.
As we have pointed out in Remark \ref{remark 3.8.01}
condition \eqref{4.9.4} implies \eqref{07.8.12.3}.
 We write the latter at $x=0$,
substitute $\varphi(x/h)$ in place of $\varphi$
and let $h\downarrow0$. Then by just comparing the powers
of $h$ in different terms we obtain that
for a constant $N$ and any $i=1,...,d$ and $\varphi$
  we have
$$
\sum_{\lambda\in\Lambda_{1}}(D_{i}q_{\lambda}(0))
\varphi(\lambda)\leq N\big(\sum_{\lambda,\nu\in\Lambda_{1}}
|\delta_{1,\lambda}
\delta_{1,\nu}\varphi(0)|^{2}\big)^{1/2}.
$$
We see that the linear function (of $\varphi$) on the left
provides a supporting plane at the origin
for the convex function on the right.
Consequently, there are some constants $q_{\lambda,
\nu}$ such that for all $\varphi$
$$
\sum_{\lambda\in\Lambda_{1}}(D_{i}q_{\lambda}(0))
\varphi(\lambda)=
\sum_{\lambda,\nu\in\Lambda_{1}}q_{\lambda,\nu}
\delta_{1,\lambda}
\delta_{1,\nu}\varphi(0)
$$
\begin{equation}
                                               \label{3.25.2}
=\sum_{\lambda,\nu\in\Lambda_{1}}q_{\lambda,\nu}
[\varphi(\lambda+\nu)+
\varphi(0)-\varphi(\lambda)-\varphi(\nu)].
\end{equation}
Without losing generality we may assume that
$$
q_{\lambda,\nu}=q_{\nu,\lambda} 
$$
and split the sum on the right in
\eqref{3.25.2} into two parts: the first part
with the summation over $\lambda,\nu$
such that $\lambda+\nu=0$ and the second part
for $\lambda+\nu\ne0$. According to assumption \eqref{3.25.4}
the terms $\varphi(\lambda+\nu)$ in the second part do 
not appear elsewhere in \eqref{3.25.2}. It follows that
$q_{\lambda,\nu}=0$ if $\lambda+\nu\ne0$, so that
$$
\sum_{\lambda\in\Lambda_{1}}(D_{i}q_{\lambda}(0))
\varphi(\lambda)
=\sum_{\lambda\in\Lambda_{1} }q_{\lambda,-\lambda}
[2
\varphi(0)-\varphi(\lambda)-\varphi(-\lambda)].
$$
Here the expression on the right is symmetric
with respect to the transformation
$\varphi(x)\to\varphi(-x)$. Thus,
$$
\sum_{\lambda\in\Lambda_{1}}(D_{i}q_{\lambda}(0))
\varphi(\lambda)=
\sum_{\lambda\in\Lambda_{1}}(D_{i}q_{ \lambda}(0))
\varphi(-\lambda).
$$
We obtained this relation at the origin. Similarly,
for any $x$
\begin{equation}
                                              \label{3.25.7}
\sum_{\lambda\in\Lambda_{1}}(D_{i}q_{\lambda}(x))
\varphi(\lambda)=
\sum_{\lambda\in\Lambda_{1}}(D_{i}q_{ \lambda}(x))
\varphi(-\lambda).
\end{equation}
Fix a $\lambda_{0}\in\Lambda_{1}$ an take   a  $\varphi$
which is $1$ as $\lambda=\lambda_{0}$ and zero otherwise.
Then \eqref{3.25.7} shows that 

(i) either $-\lambda_{0}\in\Lambda_{1}$ and then $Dq_{\lambda_{0}}(x)
=Dq_{-\lambda_{0}}(x)$,

(ii) or $-\lambda_{0}\not\in\Lambda_{1}$ and then 
$Dq_{\lambda_{0}}(x)=0$. 
 
This proves our claim.

\end{remark}

  In Remark \ref{remark 3.26.1}
we saw that \eqref{3.25.4} along with
\eqref{4.9.4} lead to
the symmetry of the operator $Q_{ \nu}$
expressed by
\eqref{3.26.2}.
 However, if $\mu+\lambda \in\Lambda_{1}$
for some $\mu,\lambda\in\Lambda_{1} $,
the symmetry of $Q_{\nu}$   may not occur.
\begin{example}
Let 
$d=1$ and $\Lambda_{1}=\{-3,-1,1,2 \}$.
Take a smooth function $f(x)$ such that
$1\leq f\leq2$ and 
set
$$
 q_{-3}=1\quad q_{-1}=q_{2}=3-f, 
\quad  q_{1}=f.
$$
Then
$$
\sum_{\lambda\in\Lambda_{1}}q_{\lambda}\lambda
=-3-(3-f)+f+2(3-f)=0.
$$
Furthermore, 
$$
h^{-1}\sum_{\lambda\in\Lambda_{1}}q'_{\lambda}\delta_{\lambda}
\varphi=h^{-2}f'(-T_{-1}+1+T_{1}-T_{2})\varphi
=f'R\varphi,
$$
where
$$
R:=h^{-2}(T_{2}-1)(T_{-1}-1)=\delta_{2}\delta_{-1}.
$$
It follows that \eqref{4.9.4} is satisfied.
Also observe that for $\mu\in\Lambda_{1}$
we have 
$$
Q_{\mu}\varphi=h^{-1}\sum_{\lambda\in\Lambda_{1}}
(\delta_{\mu}q_{\lambda})\delta_{\lambda}
\varphi=(\delta_{\mu}f)R\varphi,
$$
and
$$
Q_{\mu}T_{\mu}\varphi=Q_{\mu}(T_{\mu}-1)\varphi+
Q_{\mu}\varphi=
\sum_{\lambda\in\Lambda_{1}}
(\delta_{\mu}q_{\lambda})\delta_{\lambda}\delta_{\mu}
\varphi+Q_{\mu}\varphi.
$$
This and the fact that $q_{\lambda}\geq1$
easily imply that
condition \eqref{4.9.3} is also satisfied
with appropriate constants and operator $\cK $
in case there are terms also with $p_{\lambda}$ in $L$
and either $p_{\lambda}\geq0$ or $h$ is sufficiently small
so that $\chi_{\lambda}\geq1/2$.
\end{example}

\begin{remark}                           \label{11.17.10.07}
Note that the   argument
in the above  example 
shows that (always 
under Assumptions \ref{assumption 16.12.07.06}
and
\ref{assumption 1.26.11.06}) an operator
$$
L_h=h^{-1}\sum_{\lambda\in\Lambda_1}
q_{\lambda}\delta_{\lambda}-c
$$
satisfies Assumption \ref{assumption 11.22.11.06} 
if 
$q_{\lambda}\geq \kappa>0$  for a constant $\kappa>0$ and 
 the equality
\begin{equation}                         \label{7.15.10.07}
h^{-1}\sum_{\lambda\in\Lambda_1}
(D_iq_{\lambda})\delta_{\lambda}
=\sum_{\lambda,\mu\in\Lambda_1}q_{i\lambda\mu}
\delta_{\lambda}\delta_{\mu}, 
\quad i=1,2,\dots, d
\end{equation}                               
holds with 
  some bounded
coefficients $q_{i\lambda\mu}$. 
Therefore it would be useful to find simple conditions,  
i.e., which can be easily verified, for the characterization  
of $\Lambda_1$ and $q_{\lambda}$ satisfying 
\eqref{7.15.10.07}.  In this direction we have 
the following condition and conjecture  about
a criterion
 for 
\eqref{7.15.10.07} to hold. 
 
  We call a function $\varphi$ on
$\Lambda_1\cup\{0\}$ {\em linear\/} if
$\varphi(\nu)=\varphi(\lambda)+\varphi(\mu)$ whenever 
$\nu,\lambda,\mu\in\Lambda_1\cup\{0\}$
and $\nu=\mu+\lambda$. 

\noindent{\em Conjecture\/}. Equation \eqref{7.15.10.07}
holds with some $q_{i\lambda\mu}$ if and only if
$$
\sum_{\lambda\in\Lambda_1}(Dq_\lambda)\varphi(\lambda)=0
$$ 
holds for any $\varphi$ that is linear on $\Lambda_{1}$.

Notice that \eqref{7.15.10.07} and
the condition of the conjecture
are satisfied, for example, 
when property (S) 
holds, or if $\Lambda_1$ is the  union of
disjoint triplets 
$\{\lambda,\mu, \lambda+\mu\}$ such that 
$-Dq_{\lambda+\mu}=Dq_{\lambda}=Dq_{\mu}$ for each of them. 
\end{remark}

\begin{remark} 

Condition \eqref{07.9.25.7} is ``fool proof"
in two ways related to changes of variables.
For simplicity we only concentrate on the case
that $\tau_{0}=0$ and $\tau_{\lambda}=1$
for $\lambda\in\Lambda_{1}$.
First, one can try changing the time variable 
by   introducing the new function
$v(t,x)=u(\kappa^{-1}t,x)$, where $\kappa>0$
is a constant. This amounts to
 dividing the coefficients of
\eqref{equation} and $f$ by  $\kappa $ and
accordingly changing time $t\to\kappa^{-1} t$. 
However, as is easy to see this will not affect 
condition \eqref{07.9.25.7} and, for that matter, the value of
$F_{1}$ in Theorem \ref{theorem 4.7.2} either.

 The second way is to try to relax condition \eqref{07.9.25.7}
by changing the space
variable.
Introduce
$$
\bar{h}=\kappa h,\quad S:\psi\to S\psi(t,x)
=\psi( t,\kappa^{-1}x),
$$
$$
\bar{q}_{\lambda}=\kappa^{ 2}Sq_{\lambda},\quad
\bar{p}_{\lambda}=\kappa Sp_{\lambda},
$$
$$
\bar{L}_{\bar{h}}^{0}=\bar{h}^{-1}\sum_{\lambda
\in\Lambda_{1}}\bar{q}_{\lambda}\delta_{\bar{h},\lambda}+
 \sum_{\lambda
\in\Lambda_{1}}\bar{p}_{\lambda}\delta_{\bar{h},\lambda},
$$
$$
\bar{L}_{\bar{h},\mu}^{0}
=\bar{h}^{-1}\sum_{\lambda
\in\Lambda_{1}}(\delta_{\bar{h},\mu}
\bar{q}_{\lambda})\delta_{\bar{h},\lambda}+
  \sum_{\lambda
\in\Lambda_{1}}(\delta_{\bar{h},\mu}
\bar{p}_{\lambda})\delta_{\bar{h},\lambda}.
$$

One easily checks that
\begin{equation}
                                      \label{4.6.1}
 \delta _{\bar{h},\lambda} 
=\kappa^{-1}S\delta_{h,\lambda}S^{-1},
\quad\lambda\in\Lambda,\quad
\bar{L}_{\bar{h} }^{0}= SL_{h }S^{-1}.
\end{equation}
Owing to \eqref{4.6.1}, if $u$ satisfies 00
\eqref{equation}, then $\bar{u}=Su$ satisfies
$$
\frac{\partial\bar{u}}{\partial t}
=\bar{L}^{0}_{\bar{h}}\bar{u}-\bar{c}\bar{u}+\bar{f},
$$
where $\bar{c}= Sc$,
$\bar{f}=Sf$.
Furthermore,
 $\bar{L}_{\bar{h},\mu}^{0}=\kappa^{-1}SL_{h,\mu}S^{-1}$
and if \eqref{07.9.25.7} is satisfied, then for $\psi=S^{-1}
\varphi$
$$
2  \sum_{\lambda\in\Lambda_{1}}
(\delta_{\bar{h},\lambda}\varphi)
\bar{L}^{0}_{\bar{h},\lambda}T_{\bar{h},\lambda}\varphi
=2\kappa^{-2}S\sum_{\lambda\in\Lambda_{1}}
(\delta_{h,\lambda}\psi)
L^{0}_{h,\lambda}T_{h,\lambda}\psi
$$
$$
\leq
(1- \delta)\kappa^{-2}S\sum_{\nu\in\Lambda,\lambda \in\Lambda_{1} }
(q_{\lambda}+hp_{\lambda})|\delta_{h,\lambda}
\delta_{h,\nu}\psi|^{2}
$$
$$
+K_{1}\kappa^{-2}S\sum_{\lambda\in\Lambda_{1}}(q_{\lambda}
+hp_{\lambda})
|\delta_{h,\lambda}\psi|^{2}
+
2(1-\delta)\kappa^{-2}
S\big[c\cK_{h}\big(\sum_{\lambda\in\Lambda }
|\delta_{h,\lambda}\psi|^{2}\big) \big]  
$$
$$
=
(1- \delta) \sum_{\nu\in\Lambda,\lambda \in\Lambda_{1} }
(\bar{q}_{\lambda}+\bar{h}
\bar{p}_{\lambda})|\delta_{\bar{h},\lambda}
\delta_{\bar{h},\nu}\varphi|^{2}
$$
$$
+K_{1}\kappa^{-2} \sum_{\lambda\in\Lambda_{1}}(\bar{q}_{\lambda}
+\bar{h}\bar{p}_{\lambda})
|\delta_{\bar{h},\lambda}\varphi|^{2}
+
2(1-\delta)\bar{c} S\cK_{h}S^{-1}\big(\sum_{\lambda\in\Lambda}
|\delta_{\bar{h},\lambda}\varphi|^{2}\big),   
$$
where $S\cK_{h}S^{-1}\in\frK$. We see that this
change of coordinates did not produce any
effect on \eqref{07.9.25.7}
apart from changing $K_{1}$ and $\cK_{h}$, which is irrelevant.

\end{remark}


\begin{thebibliography}{mm}

\bibitem{DK2} Hongjie Dong and 
N.V. Krylov,
 {\em
On the rate of convergence
of finite-difference approximations for 
Bellman equations with
constant coefficients\/},  Algebra i Analiz, Vol. 
17 (2005), No. 2, 108-132;
St.~Petersburg
Math.~J, Vol. 17 (2006), No. 2, 295-313.

\bibitem{DK} Hongjie Dong and 
N.V. Krylov,
{\em On the rate of convergence
of finite-difference approximations for
degenerate linear parabolic equations
with $C^{1}$ and $C^{2}$ coefficients\/},  
Electron. J. Diff. Eqns.,
Vol. 2005(2005), No. 102, pp. 1-25.
http://ejde.math.txstate.edu

\bibitem{DK1} Hongjie Dong and 
N.V. Krylov,
{\em  On the rate of convergence
of finite-difference approximations for parabolic 
Bellman equations
with Lipschitz coefficients in cylindrical
domains\/},   Applied Math. and Optimization, 
Vol. 56 (2007), No. 1,
37-66.

\bibitem{Gy}I. Gy\"ongy, {\em
Lattice approximations for stochastic 
quasi-linear parabolic partial differential equations 
driven by space-time white noise. II\/},
Potential Anal., Vol. 11 (1999), No. 1, 1-37. 

\bibitem{GK}I. Gy\"ongy and N.V. Krylov, 
{\em Accelerated convergence of finite difference schemes 
for second order parabolic and elliptic PDEs\/},
in preparation.

\bibitem{Kr07} N.V. Krylov,
{\em A priori estimates
of smoothness of solutions to  difference  Bellman's equations 
  with linear and quasilinear operators\/}, 
Math. Comp., Vol. 76 (2007), 669-698. 

\bibitem{Kr08} N.V. Krylov, {\em On
factorizations of smooth nonnegative matrix-values
functions and on
smooth functions with values in polyhedra\/},
submitted to Applied Math. Optimiz.
  http://arxiv.org/pdf/0706.0192

\bibitem{Li} W. Littman,  {\em
R\'esolution du probl\`eme 
de Dirichlet par la m\'ethode des diff\'erences finies\/},
C. R. Acad. Sci. Paris, Vol. 247 (1958), 2270-2272.

\bibitem{Yo}Hyek Yoo, {\em
Semi-discretization of stochastic partial 
differential equations on $R\sp 1$ by a 
finite-difference method\/},
Math. Comp., Vol. 69 (2000), No. 230, 653-666. 

\end{thebibliography}
\end{document}